\documentclass{article}
\usepackage[top=1in, bottom=1in, left=1in, right=1in]{geometry}

\usepackage{amsfonts}
\usepackage{amsmath}
\usepackage{amssymb}
\usepackage{bm}
\usepackage{caption}
\usepackage{graphicx}
\usepackage{algorithm}
\usepackage{algpseudocode}
\usepackage{multirow}
\usepackage{bbm}
\usepackage{xcolor}
\usepackage{booktabs}
\usepackage{yfonts}
\usepackage{enumitem}
\usepackage{adjustbox}
\usepackage{mathtools}
\usepackage{float}
\usepackage{subfig}
\usepackage{verbatim}

\newcommand{\Be}{\bm{e}}

\newcommand{\Bp}{\bm{p}}

\newcommand{\R}{\ifmmode\mathbb{R}\else$\mathbb{R}$\fi}
\newcommand{\C}{\ifmmode\mathbb{C}\else$\mathbb{C}$\fi}
\newcommand{\N}{\ifmmode\mathbb{N}\else$\mathbb{N}$\fi}
\newcommand{\Q}{\ifmmode\mathbb{Q}\else$\mathbb{Q}$\fi}
\newcommand{\Z}{\ifmmode\mathbb{Z}\else$\mathbb{Z}$\fi}

\algnewcommand{\LeftComment}[1]{\Statex \(\triangleright\) #1}
\usepackage{hyperref}

\providecommand{\keywords}[1]
{
  \small	
  \textbf{\textit{Keywords---}} #1
}

\title{Multi-Scale Finite Expression Method for PDEs with Oscillatory Solutions on Complex Domains}

\author{
Gareth Hardwick{\thanks{Department of Mathematics, Purdue University}}, and Haizhao Yang{\thanks{Department of Mathematics and Department of Computer Science, University of Maryland College Park}}
}

\begin{document}
\maketitle
\noindent
\begin{abstract}
Solving partial differential equations (PDEs) with highly oscillatory solutions on complex domains remains a challenging and important problem. High-frequency oscillations and intricate geometries often result in prohibitively expensive representations for traditional numerical methods and lead to difficult optimization landscapes for machine learning–based approaches. In this work, we introduce an enhanced Finite Expression Method (FEX) designed to address these challenges with improved accuracy, interpretability, and computational efficiency. The proposed framework incorporates three key innovations: a symbolic spectral composition module that enables FEX to learn and represent multiscale oscillatory behavior; a redesigned linear input layer that significantly expands the expressivity of the model; and an eigenvalue formulation that extends FEX to a new class of problems involving eigenvalue PDEs. Through extensive numerical experiments, we demonstrate that FEX accurately resolves oscillatory PDEs on domains containing multiple holes of varying shapes and sizes. Compared with existing neural network–based solvers, FEX achieves substantially higher accuracy while yielding interpretable, closed-form solutions that expose the underlying structure of the problem. These advantages—often absent in conventional finite element, finite difference, and black-box neural approaches—highlight FEX as a powerful and transparent framework for solving complex PDEs.
\end{abstract}


\section{Introduction}
Partial differential equations (PDEs) serve as the mathematical backbone for modeling a wide range of phenomena, from fluid flow and electromagnetism to the behavior of financial derivatives and more. In many applications, the governing PDEs admit solutions that exhibit rapid spatial oscillations or are defined over geometrically intricate domains with nontrivial topologies, such as regions with multiple holes. Traditional approaches such as finite difference and finite element methods remain highly effective in structured, low-frequency settings, but their performance degrades as the solution becomes increasingly oscillatory or the domain geometry becomes irregular. Furthermore, the mesh size in traditional methods grows exponentially with the number of dimensions in PDEs.  This is the ``curse of dimensionality" -- in its most basic form it can be seen when estimating a general function to a given degree of accuracy - as the number of input variables increases, the computational complexity required increases exponentially.  

Deep learning has proven a powerful tool to lessen the curse of dimensionality in approximation theory \cite{Shen_2021,Shen_2021_3layers,shen2022deepnetworkapproximationachieving,maiti2025optimalneuralnetworkapproximation}, especially if a high-dimensional problem admits low-complexity structures \cite{JMLR:v25:22-0719,chen2023deep}, i.e., the number parameters required for NNs to approximation general high-dimensional functions is not exponentially in the problem dimension. This advantage motivated large advancements in the use of neural networks (NNs) to solve differential equations~\cite{KHOO_LU_YING_2021,Li_2022, liang2022stiffness,liang2024solving,lu2021machinelearningellipticpdes}. These solvers include the Deep Ritz \cite{e2017deepritzmethoddeep, pmlr-v134-lu21a}, Deep Nitsche \cite{Ming_2021}, Deep Galerkin \cite{Sirignano_2018}, and physics-informed neural network (PINN) \cite{Karniadakis2021, raissi} with substantial theoretical analysis \cite{LuoYang2024-HNA,JiaoLaiWangYangYang2023-DGMW,doi:10.1142/S021953052350015X,MishraMolinaro2023_IMAJNA_forwardPINNs,LiuHuangProtopapas2023_ResidualBasedErrorBound_PINNs}. Others have found great success reformulating the PDEs as backward stochastic differential equations \cite{ E2017, doi:10.1073/pnas.1718942115}, using NNs to approximate the gradient of the solution. These approaches offer notable advantages in terms of flexibility and scalability, and have demonstrated success in high-dimensional or irregular domains. 

However, NNs still face other challenges in solving PDEs. For examples, NNs suffer from the spectral bias in training, i.e., NNs tend to learn low-frequency components of the solution first and struggle to learn high-frequency components \cite{cai2019phaseshiftdeepneural,pmlr-v97-rahaman19a,Zhi_Qin_John_Xu_2020}. Therefore, NNs often struggle to represent highly oscillatory functions efficiently without extensive overparameterization \cite{chizat} or special activation functions \cite{liang2024reproducing, Ziqi_Liu_2020}, which can lead to poor generalization. Additionally, the black-box nature of neural networks poses difficulties when trying understand or gain intuition from the learned solutions. Finally, in a general setting, the NN training complexity in optimization (e.g., the number of iterations) \cite{na2025cursedimensionalityneuralnetwork,9321497} grows exponentially with the number of dimensions in PDEs.  This is the curse of dimensionality in the actual computational cost, though the number of parameters has no curse. 

To address these limitations, we propose a new Finite Expression (FEX) method to solve PDEs characterized by high-frequency solutions with complex geometries. In engineering applications such as acoustics, electromagnetics, or modeling structural vibrations, PDEs often exhibit highly oscillatory solutions. Standard numerical methods struggle in this regime due to mesh resolution constraints, while purely data-driven NN solvers also suffer for several challenges above. Our FEX-based method is designed to lessen these limitations as a new alternative choice.  FEX is a symbolic, mesh-free technique that constructs closed-form approximations to PDE solutions using a finite number of mathematical operators arranged into an expression tree \cite{hardwick2024solvinghighdimensionalpartialintegral,liang, song}. A key strength of the FEX framework lies in its ability to produce interpretable, compact expressions while maintaining a high degree of accuracy. Prior work has demonstrated FEX’s effectiveness in solving high-dimensional PDEs and committor problems \cite{liang, song2024finiteexpressionmethodlearning}, and the recent FEX-PG algorithm introduced a structured optimization framework for solving partial integro-differential equations with machine-level precision \cite{hardwick2024solvinghighdimensionalpartialintegral}.

Building upon the FEX-PG foundation \cite{hardwick2024solvinghighdimensionalpartialintegral}, the present work introduces a new multiscale FEX method to solve PDEs with highly oscillatory solutions on complex domains (e.g., domains containing multiple holes of varying sizes, shapes, and numbers) and eigenvalue problems. Our newly proposed FEX method introduces three key innovations:
\begin{enumerate}
    \item  Symbolic spectral composition that enables the representation and discovery of solutions containing high-frequency components.
    \item An improved linear input layer that significantly boosts FEX’s expressivity, allowing it to accurately approximate solutions that involve products of many terms.
    \item A new FEX formulation for solving eigenvalue problems, thereby broadening the scope of PDE problems solvable by FEX.
\end{enumerate}
We demonstrate how these new features allow FEX to robustly solve benchmark problems involving Helmholtz-type and Laplace-type equations with high-frequency solutions, as well as problems posed on domains with complex topological structure. Moreover, FEX consistently produces explicit symbolic expressions for the solution, providing both high accuracy and interpretability — qualities often lacking in black-box neural network approaches or discretization-based methods such as finite elements. We benchmark our method against recent neural network-based PDE solvers, highlighting its superior performance in accuracy and interpretability.  By combining symbolic learning, combinatorial optimization, and tailored architectural modifications, FEX offers a principled and powerful alternative to black-box solvers, with practical advantages in both performance and insight.

\section{Preliminaries}
This paper aims to develop a novel multiscale FEX method to solve PDEs with highly oscillatory solutions, particularly on domains with complex geometries. In addition, a new FEX formulation is introduced for eigenvalue problems, thereby extending the applicability of FEX. We begin by briefly reviewing the relevant PDEs we will solve, and by outlining the approaches in \cite{cai2023deepmartnetmartingalebased} and \cite{Ziqi_Liu_2020}, as our results will be compared with these.  The final part of this section introduces the ``vanilla" finite expression method, which is then adapted to solve these new problems.

\subsection{Partial Differential Equations}
\label{sec:pde}  
In this section we briefly introduce the selection of PDEs solved by our new FEX method inspired by \cite{hardwick2024solvinghighdimensionalpartialintegral} and \cite{liang} and introduce the functional used to evaluate the accuracy of our solution.  We begin with the Poisson equation:
\begin{equation*}
    -\Delta u(\textbf{x}) = f(\textbf{x}).
\end{equation*}
Here $\textbf{x}\in \Omega = [-1,1]^d$ and we impose Dirichlet boundary conditions
\begin{equation*}
    u(\textbf{x})|_{\partial \Omega} = g(\textbf{x}).
\label{eqn:dirichlet}
\end{equation*}
To further complicate the Poisson equation, non-linearity may be added as 
\begin{equation*}
    -\Delta u(\textbf{x}) + G(u) = f(\textbf{x}).
\label{eqn:non_lin_poisson}
\end{equation*}
The solution of the above PDE takes many different forms depending on $f(\textbf{x})$ and $G(u)$, different cases of which will be explored in the section of numerical results. The final equation solved is the Laplace Eigenvalue problem:
\begin{equation*}
    -\Delta u(\textbf{x}) = \lambda u(\textbf{x})
\label{eqn:eigenvalue}
\end{equation*}
with $\lambda$ as an unknown eigenvalue and $u(\textbf{x})$ as the unknown eigenfunction. Once again, $\textbf{x}\in \Omega = [-1,1]^d$ and we impose a boundary condition of $u|_{\partial \Omega} = 0$.  To avoid trivial solutions, we specify that $\lambda$ and $u$ are non-zero. 

To apply FEX to solve these PDEs, we propose a functional used to evaluate a candidate solution. This functional, for example, can consist of a least squares loss which combines both the loss on the domain and boundary.  For example, suppose there is a single PDE to be solved. Let $LHS(u)$ and $RHS(u)$ denote functionals representing the left and right hand sides of the PDE being solved, and define a new functional $\mathcal{D}(u) := LHS(u) - RHS(u)$. The loss on the domain is then defined as $\|\mathcal{D}(u)\|^2_{L_2(\Omega)}$. To enforce the boundary condition, the loss on the boundary is similarly defined as $\|u(\textbf{x}) - g(\textbf{x})\|^2_{L_2(\partial\Omega)}$. Thus, the loss functional $\mathcal{L}$ is given by the following:
\begin{equation*}
    \mathcal{L}(u)= \|\mathcal{D}(u)\|^2_{L_2(\Omega)} + \|u(\textbf{x})-g(\textbf{x})\|^2_{L_2(\partial\Omega)}.
\end{equation*}
This functional can be approximated using $N$ random points $(x_i)$ within the domain, where $x_i \in \Omega$, and $M$ points $(x_j)$ on the boundary with $x_j \in \partial\Omega$. Hence, we arrive at the discretized loss functional used in FEX:
\begin{align}
    \mathcal{L}(u) \approx \frac{1}{N}\sum_{i=1}^N|\mathcal{D}(\tilde{u}(x_i))|^2+ \frac{1}{M}\sum_{j=1}^M | \tilde{u}(x_j) - g(x_j)|^2.
    \label{eqn:leastsquare}
\end{align}
Clearly, $\mathcal{D}$ and $g$ are problem dependent. Once they are defined, $\mathcal{L}$ can be employed with FEX as introduced in Section~\ref{sec:fex}.  Note that an additional term is added to $\mathcal{L}$ in Section~\ref{sec:eigen} to address the challenges of eigenvalue problems in particular, but in all other cases, the loss in \eqref{eqn:leastsquare} is used without modification. The loss function in \eqref{eqn:leastsquare} is just a simple example based on the least squares idea for illustration purposes in preparation for presenting the FEX algorithm. There will be more advanced loss functions to be introduced later in Sections~\ref{sec:MBM} and \ref{sec:MNA}.

\subsection{Finite Expression Method}
\label{sec:fex}
\begin{figure}
\begin{center}
\includegraphics[scale = .08]{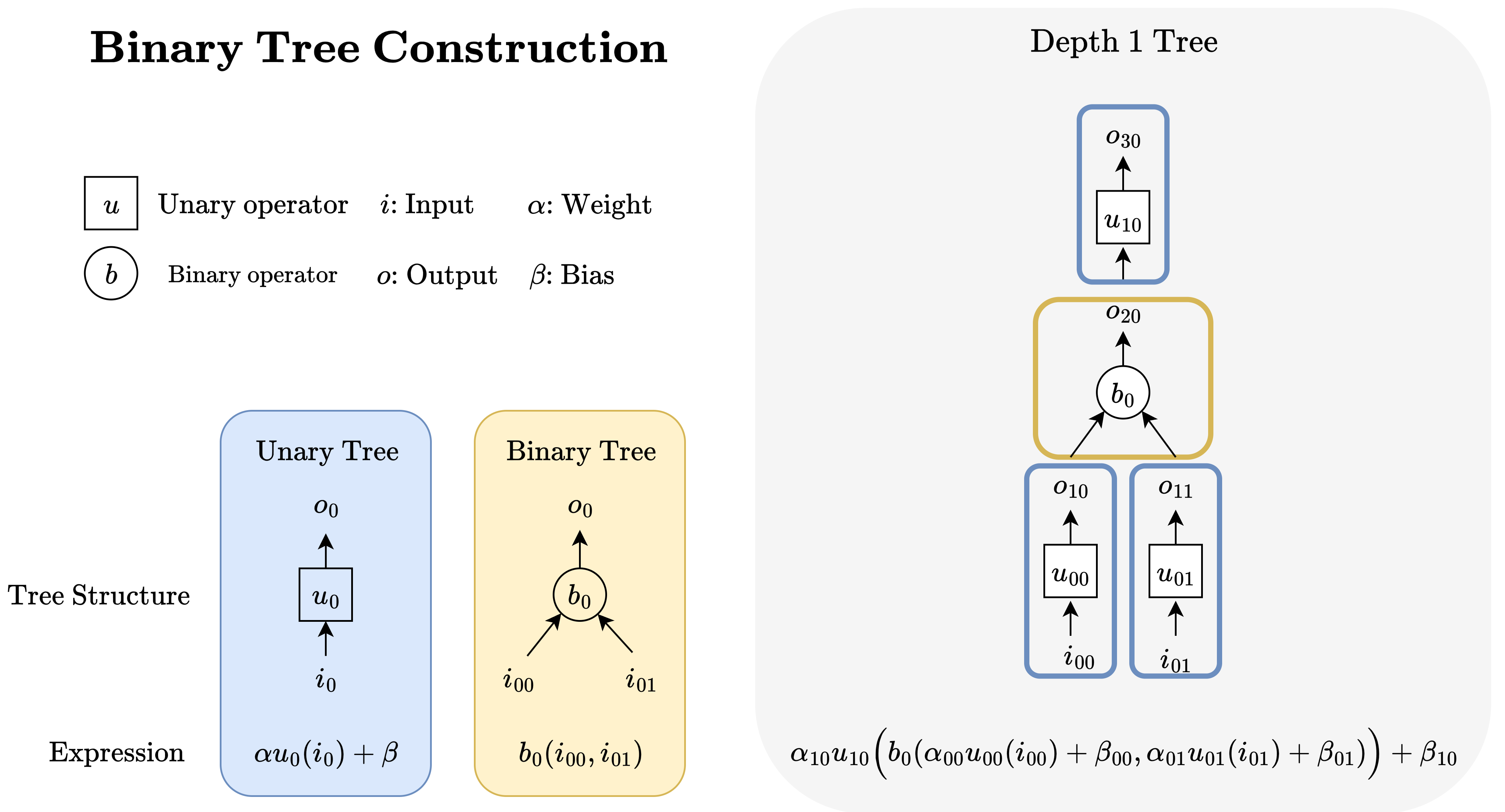}
\caption{Construction of expressions using binary trees.  Each node is either a binary or unary operator.  Beginning with the basic unary and binary trees, mathematical expressions can be built by performing computation recursively. Each tree node is either a binary operator or a unary operator that takes value from the corresponding binary or unary set. The binary set can be $\mathbb{B}:=\{+,-,\times,\cdots\}$. The unary set can be $\mathbb{U}:=\{\sin,\exp, \log, \text{Id}, (\cdot)^2, \int\cdot\text{d} x_i, \frac{\partial\cdot}{\partial x_i}, \cdots\}$, which contains elementary functions (e.g., polynomial and trigonometric function), and even integration or differentiation operators. Here ``Id'' denotes the identity map. Note that if an integration or a derivative is used in the expression, the operator can be applied using a numerical method.}
\label{fig:tree}
\end{center}
\end{figure}

FEX seeks a solution to a PDE in the function space of mathematical expressions composed of a finite number of operators. In the FEX presented, a finite expression is represented by a binary tree $\mathcal{T}$, as shown in Figure~\ref{fig:tree}. Each node in the tree is assigned a value from a set of operators, forming an operator sequence $\Be$. Each unary operator is associated with trainable weight and bias parameters, denoted by $\bm{\theta}$. Thus, a finite expression can be represented as $u(\bm{x}; \mathcal{T}, \mathbf{e}, \bm{\theta})$. The goal is to identify the mathematical expression by minimizing the functional $\mathcal{L}$, as defined in \ref{sec:pde}, where the minimizer of $\mathcal{L}$ corresponds to the solution of the PDE. Specifically, the resulting combinatorial optimization (CO) problem is:
\begin{equation*}
\min \{\mathcal{L}(u(\cdot; \mathcal{T}, \Be, \bm{\theta}))|\Be, \bm{\theta}\}.
\label{eqn:obj}
\end{equation*}

\begin{figure}
\begin{center}
\includegraphics[scale = .07]{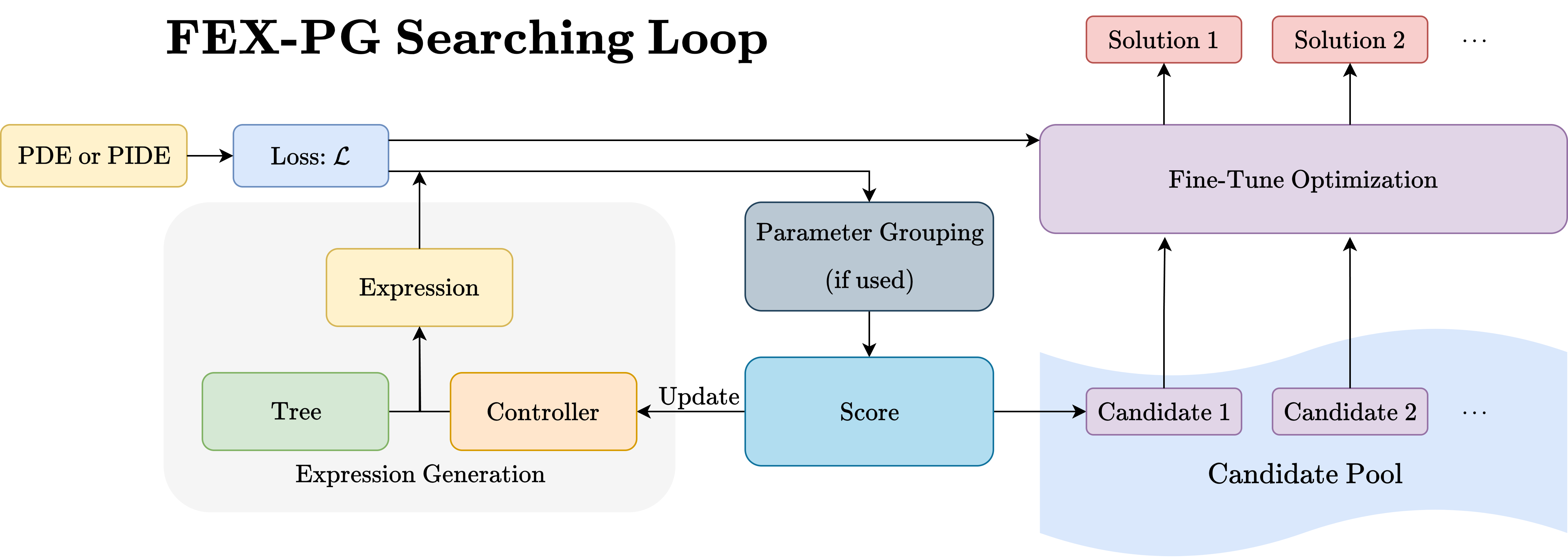}
\caption{A flowchart outlining the FEX-PG algorithm: The search loop consists of four key components: score computation, operator sequence generation, controller updates, and candidate optimization.}
\label{fig:FEXLoop}
\end{center}
\end{figure}

\begin{figure}
\centering
\includegraphics[scale = .08]{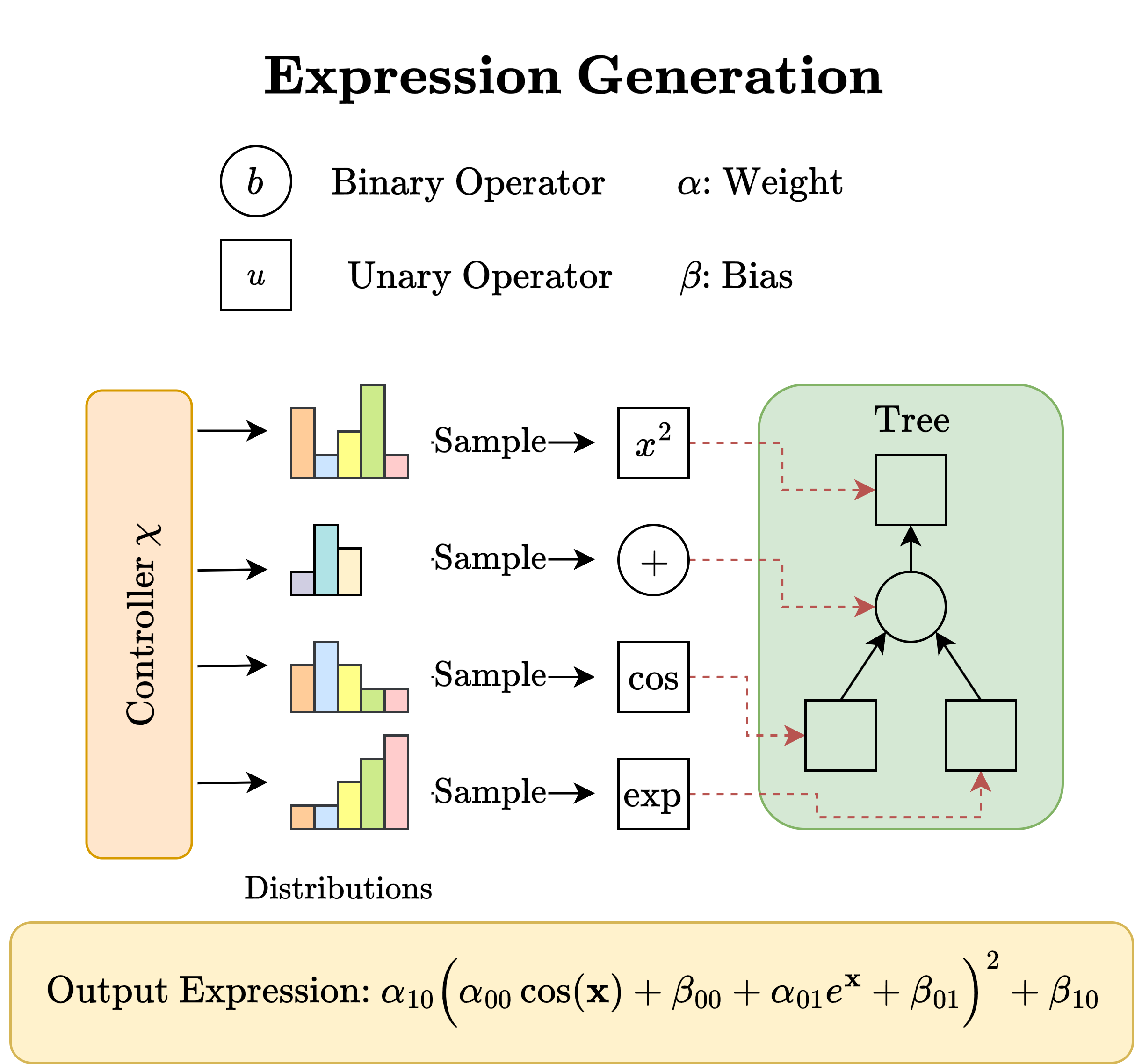}
\caption{A detailed illustration of the expression generation block from the algorithm flowchart in Figure~\ref{fig:FEXLoop}.}
\label{fig:ExpressionGen}
\end{figure}

In FEX, to address this CO, a search loop (see Figure~\ref{fig:FEXLoop}) based on reinforcement learning is employed to identify effective operators $\textbf{e}$ that can potentially recover the true solution when selected in the expression. In FEX, the search loop consists of four main components:

\begin{enumerate}
\item \textbf{Score computation (i.e., the reward in reinforcement learning)}: To efficiently evaluate the score of the operator sequence $\Be$, a mixed order optimization algorithm is used. A higher score indicates a greater likelihood that the given expression can be fine-tuned to reveal the true solution.  The score of $\Be$, denoted as $S(\Be)$, is defined on the interval $[0,1]$ by:
\begin{equation*}
S(\Be) := \big(1+L(\Be)\big)^{-1},
\label{eqn:orgscore}
\end{equation*}
where $L(\Be) := \min \{\mathcal{L}(u(\cdot; \mathcal{T}, \Be, \bm{\theta}))|\bm{\theta}\}$. As $L(\Be)$ approaches 0, the expression represented by $\Be$ approaches the true solution, causing the score $S(\Be)$ to approach 1. Conversely, as $L(\Be)$ increases, $S(\Be)$ approaches 0.
Finding the global minimizer of $\mathcal{L}(u(\cdot; \mathcal{T}, \Be, \bm{\theta}))$ with respect to $\bm{\theta}$ is computationally expensive and challenging. To speed up the computation of $S(\Be)$, rather than conducting an exhaustive search for a global minimizer using many iterations of a standard optimizer like Adam, FEX employs a combination of first-order and second-order optimization algorithms.  To begin, a first-order algorithm is employed for $T_1$ iterations to obtain a well-informed initial estimate. This well-informed initial estimate is needed so the first order method can be followed by a second-order algorithm (such as BFGS~\cite{fletcher2013practical}) for an additional $T_2$ iterations to further refine the solution.  Second order methods such as this can be highly sensitive to initial conditions, hence why the first $T_1$ iterations of the first order optimizer are critical for the success of this ``coarse-tune" process.
Let $\bm{\theta}_0^{\Be}$ denote the initial parameter set, and $\bm{\theta}_{T_1+T_2}^{\Be}$ represent the parameter set after completing $T_1+T_2$ iterations of this two-stage optimization process. The result $\bm{\theta}_{T_1+T_2}^{\Be}$ serves as an approximation of $\arg \min_{\bm{\theta}} \mathcal{L}(u(\cdot; \mathcal{T}, \Be, \bm{\theta}))$.
Then, $S(\Be)$ is estimated by:
\begin{align}
S(\Be) \approx \big(1+\mathcal{L} (u(\cdot; \mathcal{T}, \Be, \bm{\theta}_{T_1+T_2}^{\Be}))\big)^{-1}.
\label{eqn:score}
\end{align}
In FEX-PG this coarse-tune process is itself followed up by the parameter grouping process, as detailed in \cite{hardwick2024solvinghighdimensionalpartialintegral}.  The PG process serves to further refine the score of the top sequences from each batch, and sits neatly in the FEX algorithm as seen in Figure~\ref{fig:FEXLoop}.
 \item \textbf{Operator sequence generation (i.e., taking actions in RL)}:  The goal of the controller is to generate high-scoring operator sequences during the search process (see Figure~\ref{fig:ExpressionGen}). We denote the controller as $\bm{\chi}_\Phi$, where $\Phi$ represents its model parameters. Throughout the search, $\Phi$ is updated to increase the likelihood of producing operator sequences with high scores. The process of sampling an operator sequence $\Be$ from the controller $\bm{\chi}_\Phi$ is denoted as $\Be\sim\bm{\chi}_\Phi$.  Treating the tree node values of $\mathcal{T}$ as random variables, the controller $\bm{\chi}_\Phi$ outputs a number of probability mass functions $\Bp_\Phi^1, \Bp_\Phi^2, \cdots, \Bp_\Phi^s$ to characterize their distributions, where $s$ represents the total number of nodes of the tree. Each tree node value $e_j$ is sampled from its corresponding $\Bp_\Phi^j$ to generate an operator. The resulting operator sequence $\Be$ is then defined as $(e_1, e_2, \cdots, e_s)$.  The sequence is applied in-order to the tree structure, creating an expression that can be scored. To facilitate the exploration of potentially high-scoring sequences, an $\epsilon$-greedy strategy is used. With a probability of $\epsilon < 1$, $e_i$ is sampled from a uniform distribution over the set of operators. Conversely, with a probability of $1-\epsilon$, $e_i$ is sampled from $\Bp_\Phi^i$. A higher value of $\epsilon$ corresponds to a high probability of exploring new sequences, i.e. a less greedy searching process.

 \item \textbf{Controller update (i.e., policy optimization in RL)}: The controller is updated to increase the probability of generating better operator sequences based on the scores from each batch. To optimize the controller, the policy gradient method from RL is employed.  FEX makes use of the objective function proposed by \cite{petersen2021deep} to update the controller. This function is 
\begin{align*}
\mathcal{J}(\Phi)=\mathbb{E}_{\Be \sim \bm{\chi}_\Phi} \{S(\Be)|S(\Be)\geq S_{\nu, \Phi}\},
\end{align*}
where $S_{\nu, \Phi}$ denotes the $(1-\nu)\times 100\%$-quantile of the score distribution generated by $\bm{\chi}_\Phi$ within a given batch. The key detail here is that the objective function is focused only on the scores of the top performing sequences - it does not punish the controller for low scoring sequences in a given batch so long as a high scoring sequence is also present.  This is optimal in the setting of FEX since a single high scoring operator sequence is much more useful than a batch of mid-scoring ones - in theory there is likely to be a single best sequence - this is what FEX seeks to find.  In addition, this objective function also helps to avoid punishing exploration within the operator space, ensuring that exploration of new operator sequences can continue even in later iterations of the searching loop.

The controller parameters $\Phi$ are updated using gradient ascent and learning rate $\eta$: 
\begin{align*}
\Phi \leftarrow \Phi+\eta \nabla_\Phi\mathcal{J}(\Phi). 
\end{align*}

 \item \textbf{Candidate optimization (i.e., policy deployment)}: A pool of high scoring operator sequences, the ``candidate pool" is built and maintained during the search.  After this, the parameters $\bm{\theta}$ of each candidate $\Be$ in the pool are optimized to approximate the PDE solution.
The use of the pool is critical due to the challenges of scoring sequences during the searching loop.  The score of an operator sequence $\Be$ is determined by optimizing a highly nonconvex function, starting from a random initial point and using a limited number of update iterations ($T_1 + t_2)$. This approach, while efficient, can easily fail to capture the true score of a sequence if the optimization process becomes trapped at a local minima.  Because of this it is entirely possible that the operator sequence that most closely approximates or even exactly matches the true solution does not achieve the highest score.  To mitigate the risk of overlooking promising operator sequences (i.e. ones that may be able to represent the true solution, but perhaps were not scored accurately), a candidate pool $\mathbb{P}$ with a fixed capacity $K$ is used. This pool is designed to store multiple high-scoring sequences of $\Be$.

The pool is implemented as so: The top scoring sequence from each batch is added to the candidate pool.  These sequences are ordered by score within the pool itself.  Once the pool is full (i.e. once $K$ sequences are stored), if an operator sequence is found with a better score than the worst performing sequence in the pool, that worst sequence is popped from the pool and the new one is appended to it.  The sequences are once again re-ordered by  score, and the process continues iteratively.  After the search loop concludes, we perform an additional optimization step for each $\Be$ in $\mathbb{P}$, referred to as ``fine-tuning" to contrast it from the mixed order ``coarse-tuning" used earlier. Specifically, the objective function $\mathcal{L}(u(\cdot; \mathcal{T}, \Be, \bm{\theta}))$ is optimized with respect to $\bm{\theta}$ using a first-order algorithm like Adam. This optimization runs for $T_3$ iterations with a small learning rate.  Note that generally $T_3 >> T_1 + T_2$ - we can afford to be less efficient computationally since only $K$ expressions are being optimized in this time consuming manner.

\end{enumerate}


\subsection{Martingale-Based Method}\label{sec:MBM}
In the seminal works \cite{cai2023deepmartnetmartingalebased,lu,zhang2021fbsdebasedneuralnetwork}, PDEs are solved using the coupled forward backward stochastic differential equation (FBSDE).  A crucial distinction is the use of the martingale property in \cite{cai2023deepmartnetmartingalebased} to construct a functional used to optimize the model.  To begin, define the elliptic operator $\mathcal{D}$ with $\mu = \mu(\textbf{x})$ and $\sigma = \sigma(\textbf{x})$ as
\begin{align}
\mathcal{D} = \mu^T \nabla + \frac{1}{2}Tr(\sigma \sigma^T \nabla \nabla^T).
\label{eqn:ellipticop}
\end{align}
The PDE studied is given by
\begin{align}
\begin{split}
    \mathcal{D}u + V(\textbf{x},u,\nabla u) &= f(\textbf{x},u), \textbf{x}\in \Omega \subset \mathbb{R}^d, \\
    \Gamma(u) &= g, \textbf{x} \in \partial \Omega,
\label{eqn:caiPDE}
\end{split}
\end{align}
where $f$ $g$, $V$, and $\Gamma$ are problem-dependent functions and operators to specify the PDE.  For example, the boundary operator $\Gamma$ can enforce a specific condition - either the Dirichlet, Neumann, or Robin boundary condition, or a decaying condition at $\infty$ if $\Omega=\mathbb{R}^d$. Using $\mu$ and $\sigma$ to characterize drift and diffusion, respectively, \textbf{X} is sampled as the following stochastic process with a generator associated with the elliptic operator $\mathcal{D}$ from \eqref{eqn:ellipticop}: 
\begin{align*}
    &d\textbf{X}_t = \mu(\textbf{X}_t)dt + \sigma(\textbf{X}_t)\cdot d\textbf{B}_t, \\
    &\textbf{X}_0 = \textbf{x}_0 \in \Omega,
\end{align*}
where $\textbf{B}_t=(B_t^1,\dots,B_t^d)^T \in \mathbb{R}^d$ is a Brownian motion.  Rather than solving a coupled FBSDE directly as in \cite{lu}, a martingale $M^u_t$ is proposed:  
\begin{align*}
    &M^u_t = \\
    &u(\textbf{X}_t) - u(\textbf{X}_0) - \int_0^t \Big(f(\textbf{X}_s, u(\textbf{X}_s))-V(\textbf{X}_s, u(\textbf{X}_s), \nabla u(\textbf{X}_s))\Big)ds  \\
    &+ \int_0^t\Big(g(\textbf{X}_s) - cu(\textbf{X}_2)\Big)L(ds) = \int_0^t \sum_{i=1}^d \sum_{j=1}^d \sigma_{ij} \frac{\partial u}{\partial x_i}(\textbf{X}_s)dB_i(s),
\end{align*}
where $L(t)$ is the local time of the reflecting diffusion process $X_t$ (for details see \cite{2023JCoPh.47611862D}).  In the case of the Dirichlet problem, the integral with respect to local time drops out, giving a simpler form that is then used for a loss calculation. By the martingale property of $M_t \equiv M^u_t$, given a filtration $\{\mathcal{F}_s\}$ from the Brownian motion, the expectation $\mathbb{E}[M_t|\mathcal{F}_s] = M_s$. So for any measurable set $A \in \mathcal{F}_s$, we have that 
\begin{equation*}
    \mathbb{E}[M_t|A] =M_s= \mathbb{E}[M_s|A] .
\end{equation*}
Given the linearity of expectation, we then have
\begin{equation*}
    \mathbb{E}[(M_t - M_s)|A] = 0.
\end{equation*}
It is this expectation that gives rise to a functional that can be minimized to solve for the solution of the PDE \eqref{eqn:caiPDE}.  Specifically, one arrives at 
\begin{align*}
    M_t - M_s &= u(\textbf{X}_t) - u(\textbf{X}_s) - \int_s^t \mathcal{D}u(\textbf{X}_z)dz \\
    &=  u(\textbf{X}_t) - u(\textbf{X}_s) - \int_s^t \Big(f(z, u(\textbf{X}_z))-V(z, u(\textbf{X}_z), \nabla u(\textbf{X}_z))\Big)dz.
\end{align*}
In particular, if we take $A=\Omega\in \mathcal{F}_s$, we have $\mathbb{E}[\textbf{M}_t - \textbf{M}_s] = 0$, i.e., the martingale has a constant expectation. Hence, the above equation can be used to characterize the accuracy of the solution $u(\textbf{x})$ by setting the left hand equal to zero, and evaluating the right hand side across sampled trajectories.  Given a partition of the time interval $[0,T]$, $0 = t_0 < t_1 < \cdots < t_i < t_{i+1} < \cdots t_N = T$, and letting the time that the trajectory $X_t$ exits the domain be $t_{D}$, we arrive at the expression used to characterize the accuracy of a solution $u$ on the domain:
\begin{align*}
    M_{t_{i+k}\wedge t_{D}} - M_{t_i \wedge t_{D}} = u({\bf X}_{t_{i+k} \wedge t_{D}}) - u({\bf X}_{t_i \wedge t_{D}} ) - \int_{t_i \wedge t_{D}}^{t_{i+k} \wedge t_{D}}\mathcal{D}u({\bf X}_z)dz.
\end{align*}
The martingale loss \cite{cai2023deepmartnetmartingalebased} is simply the square of  $M^u_{t_{i+k}\wedge t_{D}} - M^u_{t_i \wedge t_{D}}$, averaged across the $N$ time steps of a given trajectory, i.e.,
\begin{equation}\label{eqn:lmart}
\mathcal{L}_{mart}(u):= \frac{1}{N}\sum_{i=0}^{N-1}\left(M^{u}_{t_{i+k}\wedge t_{D}} - M^{u}_{t_i \wedge t_{D}}\right)^2.
\end{equation}

\subsection{The Multi-Scale Network Approach}\label{sec:MNA}

This multi-scale network approach proposed in \cite{Ziqi_Liu_2020} is tailored specifically to PDEs with high-frequency solutions.  As such, much attention is given to learning the components of the PDE solution with different frequencies.  Once again the core concepts are presented to lend intuition to the reader.  The motivating equation solved is 
\begin{align}
    - \nabla \Big(\epsilon(\textbf{x})\nabla u(\textbf{x})\Big) + \kappa (\textbf{x})u(\textbf{x}) = f(\textbf{x}), \textbf{ x} \in \Omega \subset \mathbb{R}^d,
\label{eqn:ZiqiPDE}
\end{align}
where $\epsilon(\textbf{x})$ is the dielectric constant and $\kappa(\textbf{x})$ is the inverse Debye-Huckel length of an ionic solvent.  For simplicity, the transmission conditions on the boundary are reduced to the homogeneous boundary condition 
\begin{equation*}
    u|_{\partial \Omega} = 0.
\end{equation*}
The deep Ritz method proposed in \cite{e2017deepritzmethoddeep} is applied in \cite{Ziqi_Liu_2020} to produce a variational solution $u(\textbf{x})$ of \eqref{eqn:ZiqiPDE} (with the boundary condition above) through minimizing the energy functional as a loss function
\begin{align}
    \mathcal{L}_{Ritz}(u) = \int_\Omega \frac{1}{2}\Big(\epsilon(\textbf{x})|\nabla u(\textbf{x})|^2 + \kappa(\textbf{x})u(\textbf{x})^2\Big)d\textbf{x} - \int_\Omega f(\textbf{x})u(\textbf{x})d\textbf{x}.
    \label{eqn:energyfunctional}
\end{align}
This functional is minimized to find the variational solution as
\begin{equation*}
    u = \text{argmin}_{\nu \in H_0^1(\Omega)}\mathcal{L}_{Ritz}(\nu).
\end{equation*}
Further, the authors in \cite{Ziqi_Liu_2020} incorporated a new structure and activation function into the multi-scale network model, motivated by challenges such as the F-principle \cite{Zhi_Qin_John_Xu_2020}.  Different activation functions are tested, and the one found to be the best is 
\begin{equation*}
    \phi(x) = (x - 0)^2_+ - 3(x - 1)^2_+ + 3(x - 2)^2_+ - (x - 3)^2_+,
\end{equation*}
where $x_+ = \max \{x,0\} = \text{ReLU}(x)$.  The network is composed of many sub-networks. The input of each sub-network is a constant times the input variable $x$.  These constants ranging from $1$ to $K$ serve to specialize each sub-network to find a competent of the solution in a certain frequency range from $1$ to $K$.  Each sub-network takes $\phi(nx)$ for some $n\in[1,K]$ as the activation function, higher and lower values of $n$ corresponding to sub-networks specialized in learning higher or lower frequency components of the solution.  The output of all $K$ of these sub-networks is combined into the single output of the model, which gives the predicted value of the solution $u$ at the input variable $x$.  The model is trained using the loss function in \eqref{eqn:energyfunctional}.

\section{Multi-scale FEX for Oscillatory PDEs}
\label{sec:alg}
This section introduces new designs to make FEX capable of solving oscillatory PDEs on complex domains.  We begin with the symbolic spectral composition module, which consists of a new input layer and frequency learning strategy that, when combined, allow FEX to solve equations with oscillatory solutions. This new FEX is called the multi-scale FEX in this paper. Finally, we expand FEX further by adding to it the ability to solve eigenvalue problems.  The proposed FEX algorithm is summarized in Algorithm~\ref{alg1}.





\begin{figure}[!ht]
\begin{center}
\includegraphics[scale=.08]{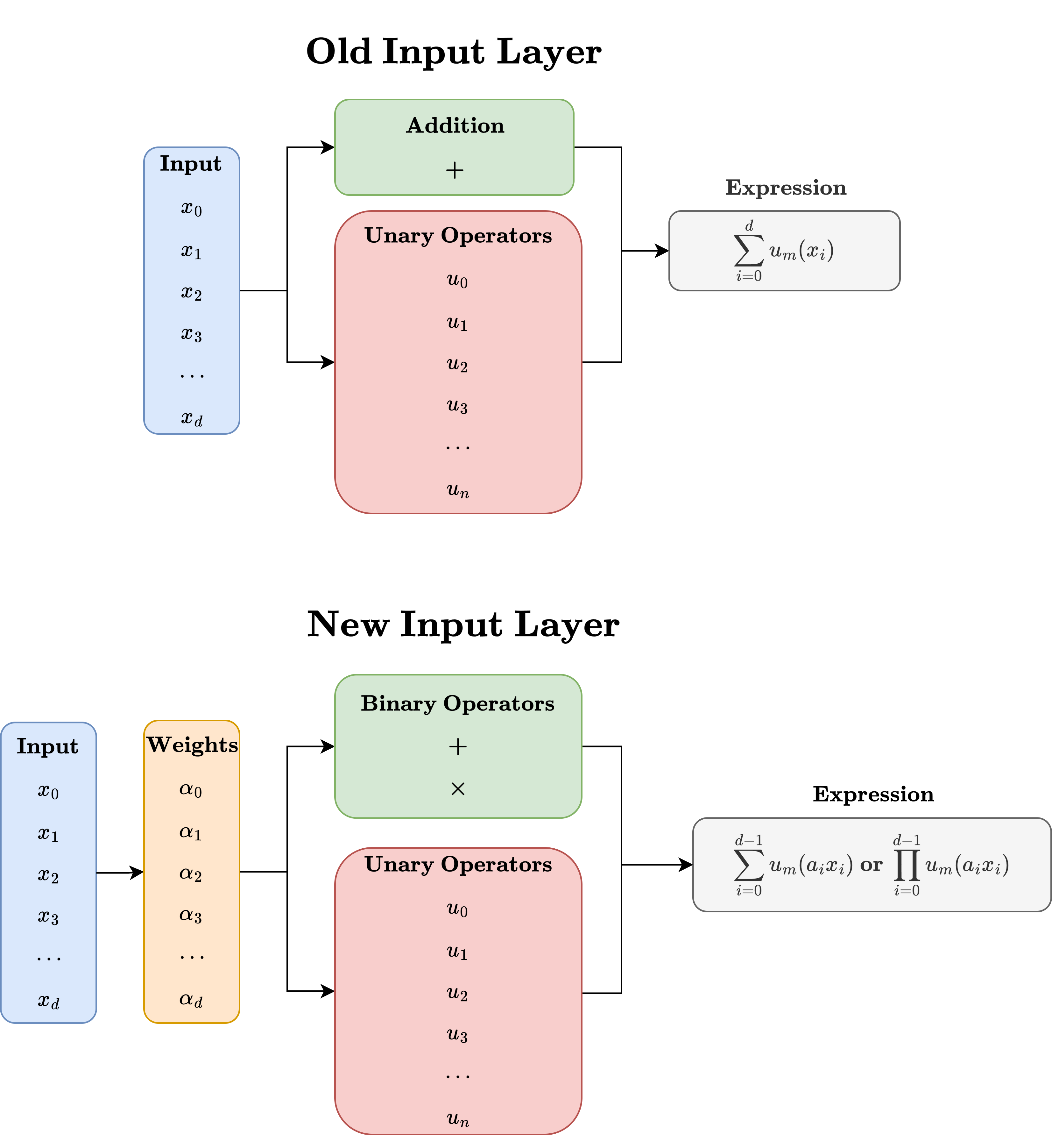}
\caption{A comparison of the old and new input layer.  Here $x\in \mathbb{R}^d$ is the input. Each component of $\textbf{x}$ is multiplied by weight $\alpha_i$. A unary operator $u_m$ is sampled from the set of unary operators in both cases, but due to the changes in the input layer structure, the resulting expressions are quite different.  The old input layer combines the terms $u_m(x_i)$ with addition always, resulting in a sum.  The new layer uses a set of binary operators (in practice, addition and multiplication), resulting in either a sum or product, depending on which operator is sampled.  This is combined with the new weights, allowing for far greater expressivity.  Both outcomes are shown to emphasize the expanded possibilities given by this reformulation.}
\label{fig:input}
\end{center}
\end{figure}

\subsection{Symbolic Spectral Composition}
\label{sec:freqlearn}
To effectively learn the frequency components of the PDE solution, two new designs are proposed to formulate the multi-scale FEX.  These designs are motivated by the main challenge of oscillatory PDEs: learning frequencies of the solution so that the correct ones can be composed together to output the true solution.  To accomplish this, a structured approach is adopted as follows.   

The first is to introduce a new input layer with two new features.  The first new feature of the input layer is an additional set of weights, implemented so that every component of the input is multiplied by a corresponding coefficient before the unary function is applied (see Figure~\ref{fig:input}).  When $u$ is chosen to be a periodic function, this allows these new coefficients (parameters $\{\alpha_i\}$ ) to determine the frequency of the periodic function.  The second new feature of the input layer is the introduction of a set of binary operators that can be sampled from.  This determines how the terms in the expression are combined into the output of the layer.  Whereas the terms were always combined using addition in the original design of FEX (creating a sum of terms, as seen at the top of Figure~\ref{fig:input}), the new layer allows for either sums or products in the new FEX proposed here.  The sampling of this binary operator occurs just the sampling of any other operator in the expression tree - the controller learns which operator to use during the expression searching loop, allowing FEX to learn the best way to combine terms (see Figure~\ref{fig:ExpressionGen} for details on expression generation and sampling).

\begin{figure}[!ht]
\begin{center}
\includegraphics[scale=.1]{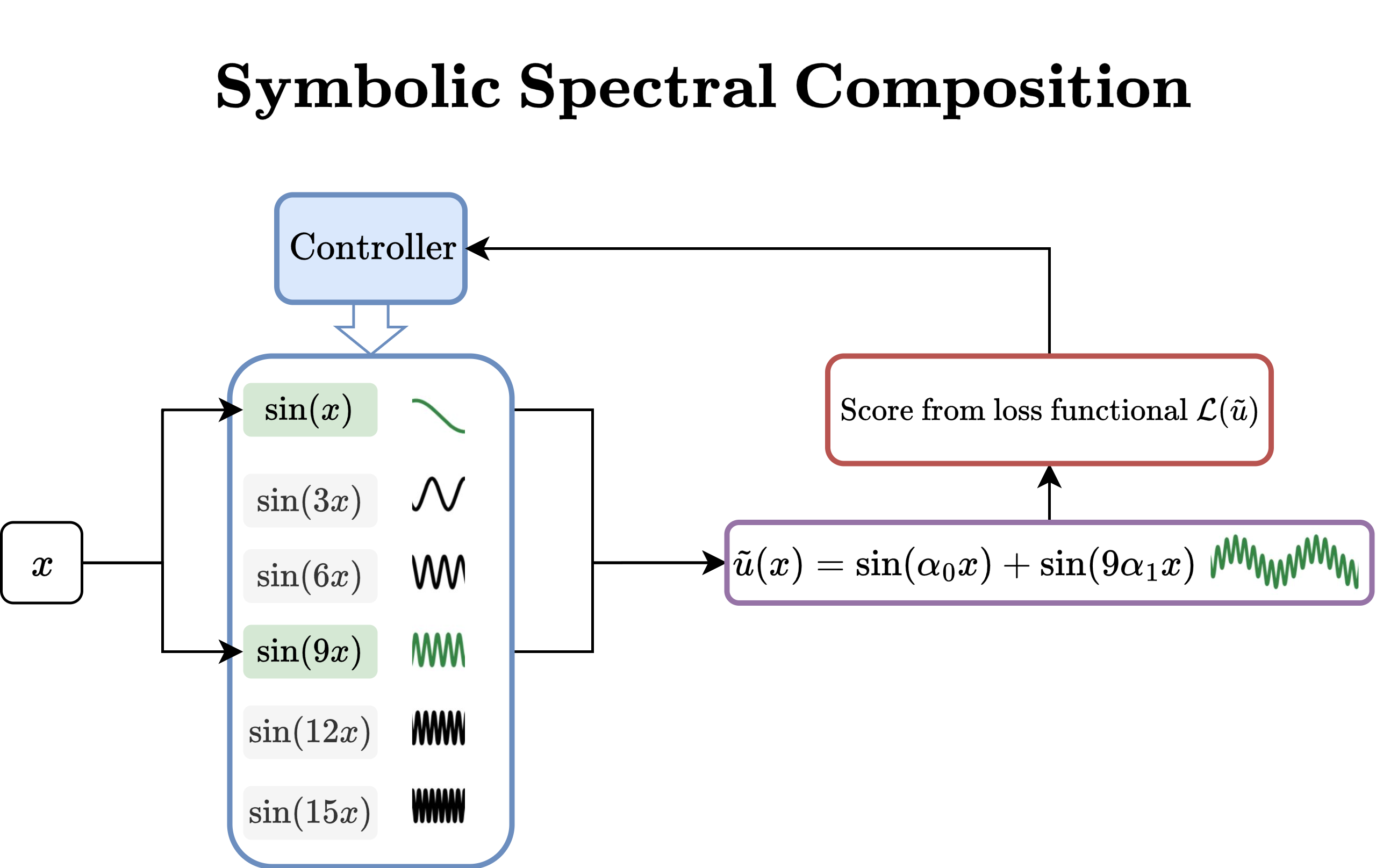}
\caption{An example of the use of the multi-scaled periodic functions, where multiple scales of sines are composed into a candidate solution $\tilde{u}(\textbf{x})$.}
\label{fig:frequencycomp}
\end{center}
\end{figure}

The second new design is to incorporate a range of scales of each periodic function into the unary operator set. Operator selection now involves two parts - which operators in general to select, and, given that a periodic operator is selected, which frequency of this operator to choose.  Following the example given in Figure~\ref{fig:frequencycomp}, the controller now selects from a set of ``base frequencies", that can then be refined by the coefficients $\{\alpha_i\}$ to minimize the loss functional $\mathcal{L}$ (in Equation \eqref{eqn:leastsquare}, \eqref{eqn:lmart}, or \eqref{eqn:energyfunctional}).  Since these coefficients are always initialized at one, the base frequencies serve as initial guesses of the true frequency components of the solution.  The search over operator sequences thus becomes a symbolic analog of constructing a spectral basis for the solution - the base frequencies  selected from the orthogonal set serve as a rough estimate of the basis, to be then adjusted and fine-tuned using the coefficients $\alpha_i$.  

The same functional $\mathcal{L}$ used to adjust the frequencies is also used to score the expressions according to Equation~\eqref{eqn:score}, which in turn is sent back to update the controller.  In this way, just as the vanilla FEX's controller that learns operators iteratively, the controller in the modified FEX also learns the base frequency composition of the solution iteratively.  Learning the base frequencies is challenging, and the score given by the loss functional is very sensitive to the base frequency used as the starting point.  To counter this, a broad spectrum of base frequencies is used to expand the set of periodic unary operators (i.e. $\sin(6x)$, $\cos(12x)$, etc.).  In practice, it was found that using our loss functional \eqref{eqn:leastsquare} with the Adam optimizer, a spacing close to $\pi$ provided an optimal balance.  While it would be tempting to include a very rich frequency spectrum - perhaps every integer frequency - this would make the operator searching process intractably long as the CO problem would be extremely large, hence the desire for a sparser set of frequencies for the basis.  

While not explored in this paper, an adjustable basis set could lend even more flexibility to the method while keeping the CO problem from becoming bloated in length.  The current implementation serves as a proof of concept - showing that the controller can serve multiple purposes at once, selecting both operators and frequencies to compose a function to fit the given problem.
 \subsection{Eigenvalue Problems}
  \label{sec:eigen}
 We next propose a new FEX for solving eigenvalue problems.  In an eigenvalue problem such as the one introduced in Section~\ref{sec:pde}, FEX must find an eigenvalue and eigenfunction pair, $(\lambda, u)$, that solves the given equation.  
 To avoid the trivial solution $u(\textbf{x}) = \lambda = 0$,  an additional normalization term is used to ensure FEX does not simply find this trivial solution. A regularization term is introduced in \cite{cai2023deepmartnetmartingalebased} as
\begin{align}
    \Big(\frac{1}{N}\sum_{i=1}^{N}|\tilde{u}(x_i)|^p -c)\Big)^2
\label{eqn:cainormalize}
\end{align}  
to augment the loss function. In \cite{cai2023deepmartnetmartingalebased}, the constants in the normalization were chosen as $c = 1$ and $p = 1$. Motivated by this design, we propose to adopt 
\begin{align*}
    \text{min}_{i \in N} \{(|\tilde{u}(x_i)|^p -c)^2\}
\end{align*}
as a new regularization term. As noted, this minimum is taken over all points in the random batch of sample points of size N.  The intuition behind the change is simple: we wish to avoid the zero solution, but it should be noted that \eqref{eqn:cainormalize} will reward solutions whose mean value is near to $c$.  This is a stronger condition than we need and causes FEX to avoid solutions that are close to zero in many places (but not trivial).  By using the minimum we punish only solutions that are near zero everywhere.  In practice, this weaker condition is seen to be sufficient for FEX to avoid the trap of the trivial solution.  The complete functional for eigenvalue problems is 
\begin{align}
    \mathcal{L}(u) \approx \frac{1}{N}\sum_{i=1}^N|\mathcal{D}(\tilde{u}(x_i))|^2+ \alpha_b \frac{1}{M}\sum_{j=1}^M | \tilde{u}(x_j) - g(x_j)|^2 + \alpha_n \text{min}_{i \in N}\{(|\tilde{u}(x_i)|^p -c)^2\},
\label{eqn:eigenloss}
\end{align}
where $\alpha_b$ and $\alpha_n$ are hyperparameters that let us weight the boundary loss and normalization loss respectively to optimize learning speed.
Note that clearly \eqref{eqn:eigenloss} is a function of both $u(\textbf{x})$ and $\lambda$ - we must solve for both simultaneously. Next, we add a learnable parameter, $\lambda$ to FEX.  Learning this eigenpair, $(u, \lambda)$, is highly sensitive to the initial guess of $\lambda$.  As such, we develop an approach to initialize our new $\lambda$ parameter.  To do this we draw inspiration from the Rayleigh Quotient \cite{haberman}.  Typically this quotient is expressed for matrices as 
\begin{align}
    R(A,x) = \frac{x^{\ast}Ax}{x^{\ast}x},
    \label{eqn:rayleigh}
\end{align}
where $A$ is a Hermitian matrix, $x$ a non-zero vector, and $\ast$ represents the conjugate transpose.  Further, if we restrict ourselves to real matrices and vectors the Hermitian condition on $A$ reduces to that of $A$ being symmetric, and the conjugate transpose $\ast$ reduces to simply the transpose $T$.  Of great importance is the equality
\begin{align*}
    R(A,x) = \frac{x^{\ast}Ax}{x^{\ast}x} = \frac{\sum_{i=1}^n \lambda_i y_i^2}{\sum_{i=1}^n y_i^2},
\end{align*}
where $(\lambda_i, v_i)$ is the $i^{th}$ eigenvalue-vector pair and $y_i = v_i^{\ast}x$ is the $i^{th}$ coordinate of $x$ in the eigenbasis.  Given this, we can easily bound this quotient above by 
\begin{align*}
    R(A,x) = \frac{x^{\ast}Ax}{x^{\ast}x} = \frac{\sum_{i=1}^n \lambda_i y_i^2}{\sum_{i=1}^n y_i^2} \leq \frac{n \lambda_{max} (v_{max}^{\ast}x)^2}{n (v_{max}^{\ast}x)^2} = \lambda_{max}.
\end{align*}
Here, $(\lambda_{max}, v_{max})$ is the eigenpair.  If $x = v_{max}$ this inequality is simply an equality and the quotient is exactly equal to the largest eigenvalue.  Equivalently from below the corresponding bound is given by the smallest eigenpair, again becoming an equality if $x = v_{min}$.  To use this concept in the context of \eqref{eqn:eigenvalue}, we write the BVP for a particular eigenpair $(u_n, \lambda_n)$, resulting in 
\begin{equation*}
    \Delta u_n + \lambda_n u_n = 0 \text{,  } u|_{\partial \Omega} = 0.
\end{equation*}
Multiply again by $u_n$ and integrate over $\textbf{x} \in \Omega$ to arrive at 
\begin{equation*}
    \int_{\Omega}u_n \Delta u_n d\textbf{x} + \lambda_n \int_{\Omega} (u_n)^2d\textbf{x} = 0.
\end{equation*}
Rearranging and using integration by parts yields
\begin{align*}
    \lambda_n \int_{\Omega} (u_n)^2d\textbf{x} &= -\int_{\Omega}u_n \Delta u_n d\textbf{x} \\
    &= \int_{\Omega} |\nabla u_n(\textbf{x})|^2d\textbf{x} - \int_{\partial \Omega} u_n(\textbf{x})\frac{\partial u_n}{\partial \nu}dS(\textbf{x}).
\end{align*}
By the boundary value assumption however, the surface integral on the right hand side must be zero, and we are left with 
\begin{equation*}
     \lambda_n \int_{\Omega} (u_n)^2d\textbf{x} =\int_{\Omega} |\nabla u_n(\textbf{x})|^2d\textbf{x}.
\end{equation*}
Letting this be greater than zero (i.e. assuming $u_n$ is not constant and therefore non-trivial since $u_n|_{\partial \Omega} = 0$), lets us rearrange to arrive at 
\begin{equation*}
    \lambda_n = \frac{\int_{\Omega} |\nabla u_n(\textbf{x})|^2d\textbf{x}}{ \int_{\Omega} (u_n)^2d\textbf{x}}.
\end{equation*}
To apply this in our case, rather than $u_n(\textbf{x})$ we use $\tilde{u}(\textbf{x})$, one of our candidate functions.  We evaluate the quotient and get an initial estimate for $\lambda$ which we will call $\tilde{\lambda}_0$.  Since we sample $N$ points on the domain $\Omega$ we then arrive at our method of initializing the eigenvalue parameter:
\begin{align}
    \tilde{\lambda}_0 = \frac{\frac{1}{N}\sum_{i=1}^{N}|\nabla \tilde{u}(\textbf{x}_i)|^2}{\frac{1}{N}\sum_{i=1}^{N}(\tilde{u}(\textbf{x}_i))^2}.
\label{eqn:lambda_init}
\end{align}
Since the derivatives of $u$ are already used to compute loss in \eqref{eqn:eigenloss}, this computation is extremely efficient and does not increase complexity.  The algorithm remains identical to algorithm \ref{alg1}, except that now the functional $\mathcal{L}$ has an additional parameter, $\lambda$, initialized as in \eqref{eqn:lambda_init}.  
\begin{algorithm}
\begin{algorithmic}[1]
\caption{Fixed-Tree FEX-PG for PDEs}
\LeftComment{Input: PDE; A tree $\mathcal{T}$; Searching loop iteration $T$; Coarse-tune iteration $T_1$ with Adam; Coarse-tune iteration $T_2$ with LBFGS; Medium-tune iteration $T_3$ with Adam; Fine-tune iteration $T_4$ with Adam; Pool size $K$; Batch size $N$; Clustering threshold $\eta$.}
\LeftComment{Output: The solution $u(\textbf{x}; \mathcal{T}, \hat{\Be}, \hat{\theta})$}
\State Initialize the agent $\chi$ for the tree $\mathcal{T}$
\State $\mathbb{P} \leftarrow \{\}$
\For{$\hbox to 1em{\thinspace\hrulefill\thinspace}$ from 1 to $T$}
    \State Sample $N$ sequences $\{\Be^{(1)}, \Be^{(2)},...,\Be^{(N)}\}$ from $\chi$
    \State Losses $\leftarrow [\text{ }]$
    \For{n from 1 to $N$}
        \State Minimize $\mathcal{L}(u(\cdot; \mathcal{T}, \Be^{(n)}, \theta^{(n)}))$ with respect to $\theta^{(n)}$ by coarse-tune with $T_1 + T_2$ iterations
        \State After $T_1 + T_2$ iterations, Losses.append($\mathcal{L}(u(\cdot; \mathcal{T}, \Be^{(n)}, \theta_{T_1 + T_2}^{(n)}))$)
    \EndFor
    \State Denote $ \Tilde{n} := \arg\min (\text{Losses})$
    \State Apply operator sequence $\Be^{(\Tilde{n})}$ to tree $\mathcal{T}$, denoted as $\mathcal{T}_{e^{(\Tilde{n})}}$
    \For{leaf in $\mathcal{T}_{e^{(\Tilde{n})}}$} \Comment{Parameter Grouping}
        \State  Apply hierarchical clustering algorithm with threshold parameter $\eta$
        \State Replace the linear layer of each leaf with the modified linear layer (see \cite{hardwick2024solvinghighdimensionalpartialintegral} for details)
    \EndFor
    \For{$\hbox to 1em{\thinspace\hrulefill\thinspace}$  from 1 to $T_3$} \Comment{Learning weights for new modified linear layers}
        \State Calculate $ \mathcal{L}(u(\cdot;\mathcal{T}_{e^{(\Tilde{n})}}, e^{(\Tilde{n})}, \theta^{(\tilde{n})}))$ using $\mathcal{T}_{e^{(\Tilde{n})}}$ and update $\theta$ with Adam
        \If{$\hbox to 1em{\thinspace\hrulefill\thinspace} = T_3$ and Losses[$\tilde{n}$] $<$ $\mathcal{L}(u(\cdot;\mathcal{T}_{e^{(\Tilde{n})}}, e^{(\Tilde{n})}, \theta_{T_3}^{(\tilde{n})}))$}
            \State $\text{Losses}[\Tilde{n}] \leftarrow  \mathcal{L}(u(\cdot;\mathcal{T}_{e^{(\Tilde{n})}}, e^{(\Tilde{n})}, \theta_{T_3}^{(\tilde{n})}))$
            \EndIf
    \EndFor
    \State Calculate rewards using Losses[:] and update $\chi$
    \For{n from 1 to $N$}        
        \If{Losses$[n] <$ any in $\mathbb{P}$}
            \State $\mathbb{P}$.append($\Be^{(n)}$)
            \State $\mathbb{P}$ pops $\Be$ with the smallest reward when overloading
            \EndIf
    \EndFor
\EndFor
\For{$\Be$ in $\mathbb{P}$} \Comment{Candidate optimization}
    \For{$\hbox to 1em{\thinspace\hrulefill\thinspace}$  from 1 to $T_4$}
        \State Minimize $\mathcal{L}(u(\cdot; \mathcal{T}, \Be, \theta))$ with respect to $\theta$ using Adam
    \EndFor
\EndFor
\State \textbf{Return} the expression with the smallest fine-tune error
\label{alg1}
\end{algorithmic}
\end{algorithm}
\section{Numerical Results}
\label{sec:results}
 This section presents results showing the effectiveness of the proposed methods for solving oscillatory PDEs and eigenvalue problems on complex domains. We begin with more simple examples, and work up to complex cases with high-frequency solutions and complicated domain geometries.  First, two examples of the Poisson-Boltzmann equation from \cite{cai2023deepmartnetmartingalebased} are solved.  These initial examples are more simple and serve to validate the usage of the multi-scale FEX for these problems involving oscillatory solutions before the method is pushed further.  Then we present a number of equations from \cite{Ziqi_Liu_2020} that test the ability of multi-scale FEX to learn solutions involving high frequencies on complex domains.  All of these problems are solved using the simple least squares loss function \eqref{eqn:leastsquare}.  Finally, we end with an eigenvalue problem from \cite{cai2023deepmartnetmartingalebased}, where we then use the modified loss function (which includes the normalization term) \eqref{eqn:eigenloss}.  Picking examples from existing literature gives an apples-to-apples comparison to demonstrate the strengths of FEX. In all of these examples, the set of binary operators used is: ``$+$", ``$-$" and ``$\times$".  The unary operators are commonly used operators: ``$0$", ``$1$", ``$x$", ``$x^2$'', ``$x^3$", ``$x^4$", ``$e^x$", ``$\sin (x)$ ",``$\cos (x)$" but now augmented with multi-scale variants ``$\sin(3x)$", ``$\sin(6x)$", ... , ``$\sin(24x)$"  and ``$\cos(3x)$", ``$\cos(6x)$", ... , ``$\cos(24x)$", which allows the learning of high-frequency solutions as discussed in Section~\ref{sec:freqlearn}. Along with these operators, a depth-two tree structure is used, which can be seen on the far right of Figure~\ref{fig:tree}.  While we primarily use absolute relative error as our metric of accuracy, we also compute the $L^2$ and relative $L^2$ errors so that our results can be directly compared to those found in \cite{cai2023deepmartnetmartingalebased} and \cite{Ziqi_Liu_2020}.

\subsection{Poisson-Boltzmann Equation}
\label{sec:1d}
Here we solve the Dirichlet boundary value problem (BVP) of the Poisson-Boltzmann equation.  These examples are found in \cite{cai2023deepmartnetmartingalebased}, to which we also compare our results. \\ \\
\textbf{Example 1. } The first BVP is given by
\begin{align} \begin{cases} 
    \Delta u(\textbf{x}) + cu(\textbf{x}) = f(\textbf{x}), & \textbf{x} \in  \Omega,\\
    u(\textbf{x}) = g(\textbf{x}), &\textbf{x} \in \partial \Omega,\\
    \end{cases}
\label{eqn:PB_ex1}
\end{align}
where $c = -1$ in our case.  As in \cite{cai2023deepmartnetmartingalebased}, the true solution is
\begin{equation*}
    u(\textbf{x}) = \sum_{i=1}^d \cos(\omega x_i) \text{ with } \omega = 2.
\end{equation*}
by choosing $g$ and $f$ appropriately. Here $g(\textbf{x}) = u(\textbf{x})$ (i.e. $g$ is the true solution), and $f(\textbf{x}) = \Delta u(\textbf{x}) + cu(\textbf{x}) = -5u(\textbf{x})$. Plugging in this example to the loss function \eqref{eqn:leastsquare}, we arrive at the functional used for this problem:
\begin{equation}
    \mathcal{L}(\tilde{u}) := \frac{1}{N}\sum_{i=1}^N| \Delta\tilde{u}(\textbf{x}_i) - \tilde{u}(\textbf{x}_i) - f(\textbf{x}_i)|^2+ \frac{1}{M}\sum_{j=1}^M | \tilde{u}(\bf{x}_j) - g(\textbf{x}_j)|^2.
\end{equation}
Note that $\tilde{u}$ refers to the candidate solution generated by FEX.  
\begin{figure}[H]
    \centering
    \subfloat[Loss During Fine-tuning]{%
        \includegraphics[width=6cm]{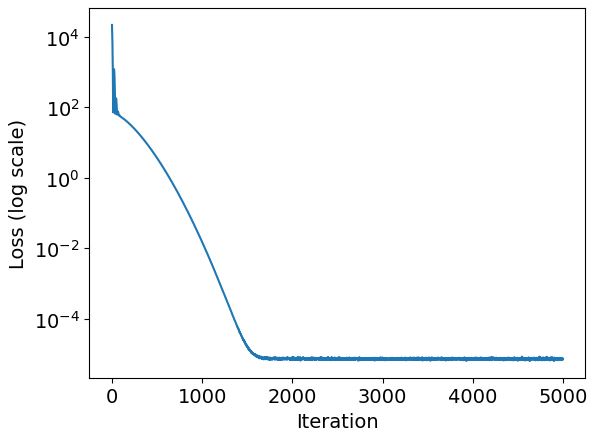}%
    }\qquad
    \subfloat[Relative Error During Fine-Tuning]{%
        \includegraphics[width=6cm]{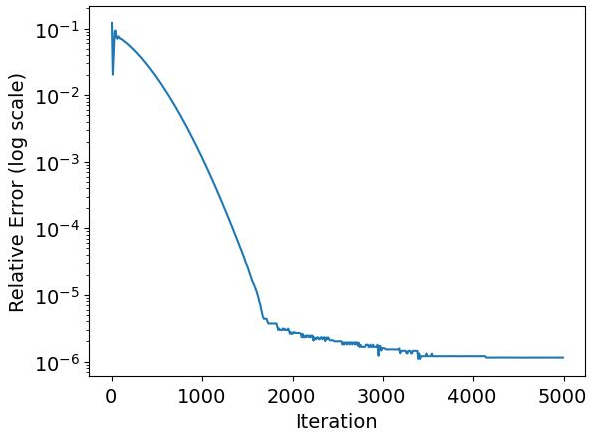}%
    }%
    \caption{Optimization profile for PDE~\eqref{eqn:PB_ex1}. 
    \textbf{(a)} Training loss of the candidate solution during fine-tuning. 
    \textbf{(b)} Absolute relative error of the candidate solution.}%
    \label{fig:5.1Test1loss}
\end{figure}
We solve the equation in \eqref{eqn:PB_ex1} on a 100-dimensional unit sphere, centered at the origin.  Note that our multi-scale FEX exhibits excellent performance in high-dimensional problems. After just $T=50$ iterations in the search loop in Algorithm \ref{alg1}, many promising mathematical expressions as candidate solutions are identified. During the fine-tuning stage, Figure~\ref{fig:5.1Test1loss} (a) shows the loss of the candidate function as it is fine-tuned over iterations, while Figure~\ref{fig:5.1Test1loss} (b) displays the corresponding absolute relative error. It is evident that after only 2000 iterations of fine-tuning, the absolute relative error drops to the order of $10^{-6}$. \\ \\
\textbf{Example 2. } The next BVP is given by
\begin{align} 
\begin{cases} 
    -\Delta u(\textbf{x}) + \sinh(u(\textbf{x})) = f(\textbf{x}), & \textbf{x} \in  \Omega,\\
    u(\textbf{x}) = g(\textbf{x}), &\textbf{x} \in \partial \Omega,\\
\end{cases}
\label{eqn:PB_ex2}
\end{align}
where $\Omega = \{\textbf{x} \in \mathbb{R}^d: ||\textbf{x}||_2 \leq 1\}$ as the unit ball.  As in \cite{cai2023deepmartnetmartingalebased}, the true solution is given as
\begin{equation*}
    u(\textbf{x}) = 2 \sum_{i=1}^d x_i^2
\end{equation*}
by choosing $f$ and $g$ appropriately. Here $g(\textbf{x}) = u(\textbf{x})$ (i.e. $g$ is the true solution) and  $f(\textbf{x}) = -\Delta u(\textbf{x}) + \sinh(u(\textbf{x})) = -4d + \sinh(u(\textbf{x}))$, where $d$ refers to the number of dimensions of $\textbf{x}$. Based on the least-square idea in \eqref{eqn:leastsquare}, our loss function is designed as:
\begin{equation}
    \mathcal{L}(\tilde{u}) := \frac{1}{N}\sum_{i=1}^N| -\Delta\tilde{u}(\textbf{x}_i) + \sinh(\tilde{u}(\textbf{x}_i)) - f(\textbf{x}_i)|^2+ \frac{1}{M}\sum_{j=1}^M | \tilde{u}(\bf{x}_j) - g(\textbf{x}_j)|^2.
\end{equation}

\begin{figure}[H]
    \centering
    \subfloat[Loss During Fine-tuning]{%
        \includegraphics[width=6cm]{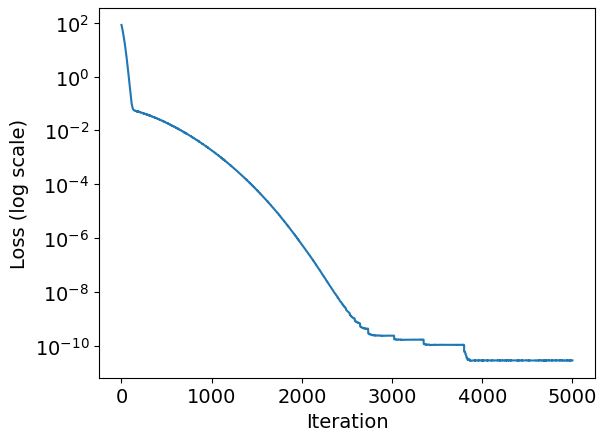}%
    }\qquad
    \subfloat[Relative Error During Fine-Tuning]{%
        \includegraphics[width=6cm]{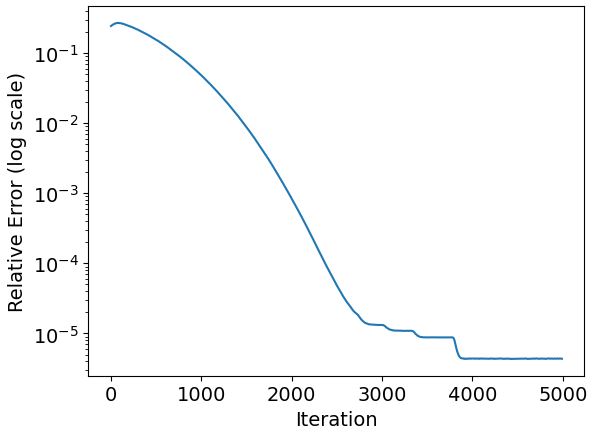}%
    }%
    \caption{Optimization profile for PDE~\eqref{eqn:PB_ex2}. 
    \textbf{(a)} Training loss of the candidate solution during fine-tuning. 
    \textbf{(b)} Absolute relative error of the candidate solution.}%
    \label{fig:5.1Test4loss}
\end{figure}

This equation was solved in a 10-dimensional unit ball.  Again, our multi-scale FEX exhibits strong performance with the absolute relative error dropping rapidly during fine-tuning. To compare with the results of \cite{cai2023deepmartnetmartingalebased} by neural networks, we also compute the average relative $L^2$ error across trials by our method. The average relative $L^2$ error by our method is 3.3e-6 in this example, comparing favorably to the method in \cite{cai2023deepmartnetmartingalebased}, where the relative $L^2$ error was around 2.5e-1. 

\subsection{Poisson Equation on a Complex 2-D Domain}
\begin{figure}[H]
    \centering
    \subfloat[Loss During Fine-tuning]{%
        \includegraphics[width=6cm]{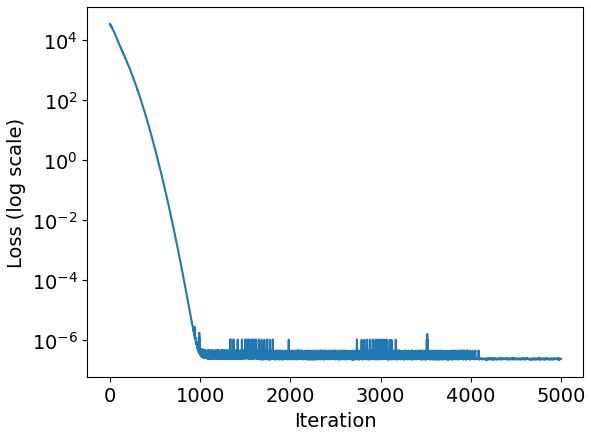}%
    }\qquad
    \subfloat[Relative Error During Fine-tuning]{%
        \includegraphics[width=6cm]{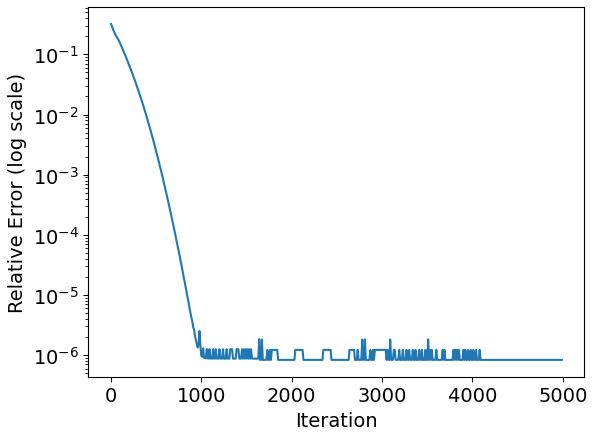}%
    }%
    \caption{Optimization profile for PDE~\eqref{eqn:Poisson6.1.3} on the first complex domain. 
    \textbf{(a)} Training loss of the candidate solution during fine-tuning. 
    \textbf{(b)} Absolute relative error of the candidate solution.}
    \label{fig:Poisson6.1.3Dom1Optim}
\end{figure}

\begin{figure}[H]
    \centering
    \subfloat[Loss During Fine-tuning]{%
        \includegraphics[width=6cm]{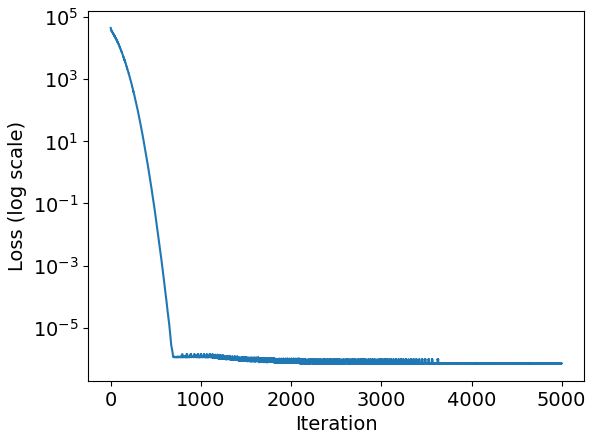}%
    }\qquad
    \subfloat[Relative Error During Fine-tuning]{%
        \includegraphics[width=6cm]{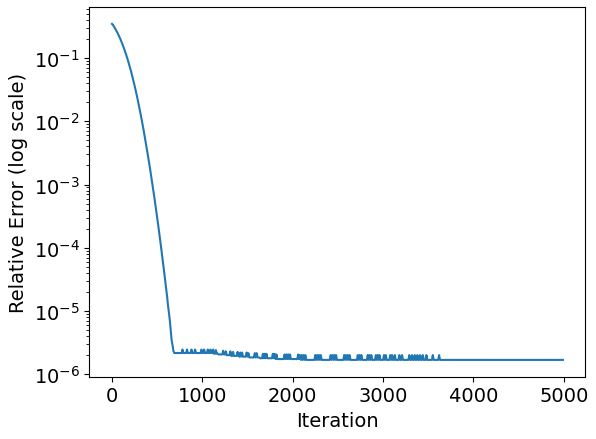}%
    }%
    \caption{Optimization profile for PDE~\eqref{eqn:Poisson6.1.3} on the second complex domain. 
    \textbf{(a)} Training loss of the candidate solution during fine-tuning. 
    \textbf{(b)} Absolute relative error of the candidate solution.}
    \label{fig:Poisson6.1.3Dom2Optim}
\end{figure}

Next, we turn our attention to examples that will validate the performance of the multi-scale FEX on complex domains, e.g., a domain featuring large holes as a test example in \cite{Ziqi_Liu_2020}.  The equation being solved is
\begin{align}
    -\Delta u(\textbf{x}) = 2\mu^2\sin(\mu x_1)\sin(\mu x_2)
    \label{eqn:Poisson6.1.3}
\end{align}
with an appropriate Dirichlet boundary condition such that the true solution is
\begin{equation*}
    u(\textbf{x}) = \sin(\mu x_1)\sin(\mu x_2),
\end{equation*}
where we let the frequency parameter $\mu$ be $7\pi$.  The PDE domain is a square with multiple holes (see Figure \ref{fig:Poisson6.1.3Dom1Optim} for an example). Based on the least squares idea in \eqref{eqn:leastsquare}, the loss functional in this example is
\begin{equation}
\mathcal{L}(\tilde{u}) := \frac{1}{N}\sum_{i=1}^N| -\Delta\tilde{u}(x_{i_1},x_{i_2}) - 2 \mu^2 \sin(\mu x_{i_1})\sin(\mu x_{i_2})|^2+ \frac{1}{M}\sum_{j=1}^M | \tilde{u}(x_{j_1},x_{j_2}) - \sin(\mu x_{j_1})\sin(\mu x_{j_2})|^2.
\end{equation}
The equation is solved on two different domains to observe the response of the multi-scale FEX to holes of different size.  The first domain has three holes centered at $-(0.5,-0.5)$, $(0.5,0.5)$, and $(0.5,-0.5)$ with radii $0.1$, $0.2$, and $0.2$, respectively (see Figure \ref{fig:Poisson6.1.3Dom1Optim}).  The second domain features four holes (see Figure \ref{fig:Poisson6.1.3Dom2Optim}).  Three of which are circles centered at $(-0.6, -0.6)$, $(0.3, -0.3)$, and $(0.6, 0.6)$ with radii $0.3$, $0.6$, and $0.3$, respectively, and the fourth one is an ellipse described by $16(x_1 + 0.5)^2 + 64(x_2 - 0.5)^2 = 1$.  All of these are exactly as in \cite{Ziqi_Liu_2020} to keep comparison fair.
\begin{figure}[H]
    \centering
    \includegraphics[width=15cm]{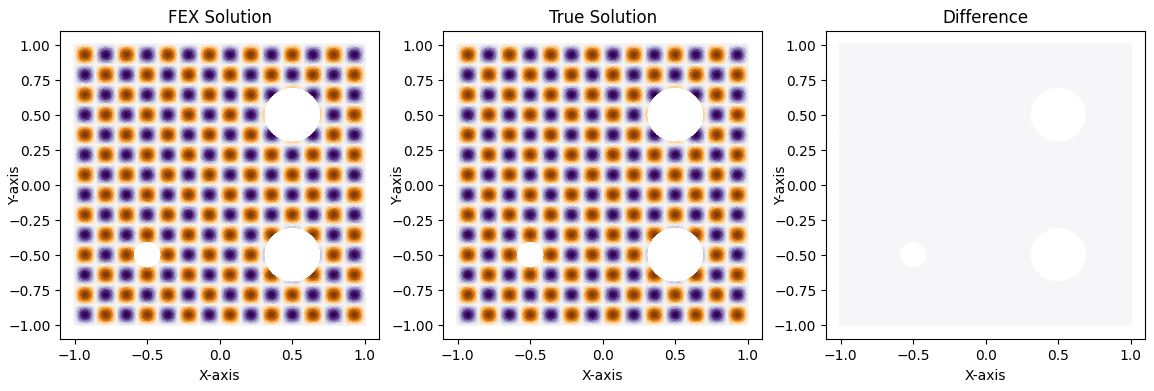}
    \caption{Comparison of FEX and exact solution on the first complex domain.  
    The rightmost figure shows the absolute difference between the true solution and the FEX solution.}
    \label{fig:Poisson6.1.3Dom1}
\end{figure}

The average absolute relative error over ten trials was $8.3 \times 10^{-7}$, with an average relative $L_2$ error of $4.9 \times 10^{-7}$. In contrast, the results of \cite{Ziqi_Liu_2020} reported a relative $L_2$ error of approximately $1 \times 10^{-2}$, highlighting the significant accuracy advantage of the multi-scale FEX. As illustrated in Figure~\ref{fig:Poisson6.1.3Dom1}, the solution produced by the multi-scale FEX is virtually indistinguishable from the true solution. A trial solution identified by the multi-scale FEX in this case was the function 
\begin{equation}
\begin{aligned}
    u(x_1, x_2) = (0.9950\sin(24(0.9162x_1))+0.0000\sin(24(0.6849x_2))+0.0000) \times \\(0.0000\sin(21(0.8286x_1))+1.0050\sin(21(1.0471x_2))+0.0000),
\end{aligned}
\end{equation}
which simplifies to
\begin{align}
    u(x_1,x_2) = 0.9999 \sin(21.9911 x_1) \sin(21.9911 x_2).
\label{eqn:FEXoutputPoisson6.1.3}
\end{align}
Note that the true solution is $u(x_1,x_2) = \sin(7\pi x_1)\sin(7\pi x_2)$.  While we have omitted the decimal expansions in the examples above for brevity, it is worth noting that the coefficients within the $\sin$ functions are accurate to six decimal places in the expansion of $7\pi$.
\begin{figure}[H]
    \centering
    \includegraphics[width=15cm]{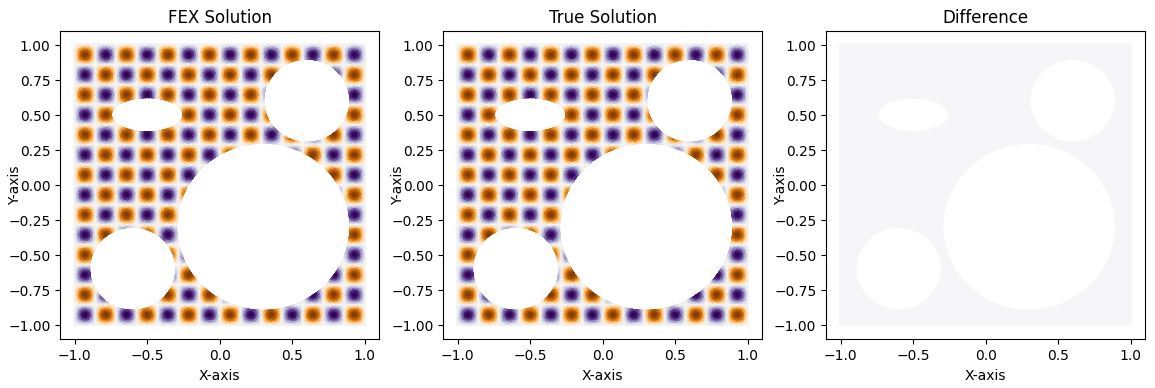}
    \caption{Comparison of FEX and exact solution on the second complex domain.  
    The rightmost figure shows the absolute difference between the true solution and the FEX solution.}
    \label{fig:Poisson6.1.3Dom2}
\end{figure}

Figure~\ref{fig:Poisson6.1.3Dom2} demonstrates that the multi-scale FEX produces a solution virtually identical to the exact one. Notably, the method shows strong robustness to the presence of holes in low-dimensional problems. The training times across the cases were nearly identical. The absolute relative error remains extremely small at $1.6 \times 10^{-6}$, and the relative $L_2$ error of the multi-scale FEX is $8.6 \times 10^{-7}$, representing a substantial improvement in accuracy compared to the results of \cite{Ziqi_Liu_2020}, which reported a relative $L_2$ error of approximately $8 \times 10^{-3}$.

\subsection{Poisson Equation on a Complex 3-D Domain}
\begin{figure}[H]
    \centering
    \includegraphics[width = 6cm]{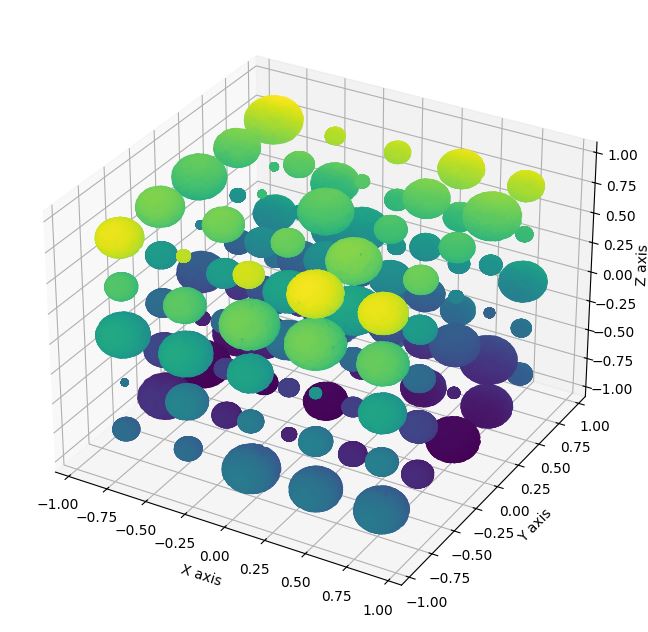}
    \caption{Cubic domain with many holes.  The coloring is used to show depth and dimensionality.}
    \label{fig:Poisson6.1.4domain}
\end{figure}

The next two examples solve Poisson equations with different true solutions corresponding to different right-hand-side functions $f(\textbf{x})$.  The PDE domain is a cube with side length $L=2$, centered on the origin, and many spherical holes inside the cube. In particular, 125 holes with random radii are placed evenly on a grid in the cube, as seen in Figure~\ref{fig:Poisson6.1.4domain}.  In the following examples, we test the maximum accuracy achievable by FEX in such domains.  To this end, we forgo a bit of speed and use double precision floats for all calculations.  The results soundly demonstrate that, in a low-dimensional setting, FEX's limit in these problems is the floating point error itself.
\noindent
\textbf{Example 1.}  The first example solved in this domain is given by
\begin{align}
    -\Delta u(\textbf{x}) = 3 \mu^2 \sin( \mu x_1)\sin(\mu x_2)\sin(\mu x_3),
    \label{eqn:Poisson6.1.4ex1}
\end{align}
with the boundary condition on the sides of the cube and surface of the spherical holes appropriately chosen such that the true solution is 
\begin{align}
    u(\textbf{x}) = \sin( \mu x_1)\sin(\mu x_2)\sin(\mu x_3).
    \label{eqn:Poisson6.1.4ex1sol}
\end{align}
Here we match \cite{Ziqi_Liu_2020} and again set $\mu = 7\pi$.  Based on the least squares idea in \eqref{eqn:leastsquare}, the loss functional in this example is
\begin{align*}
\mathcal{L}(\tilde{u}) & := \frac{1}{N}\sum_{i=1}^N| -\Delta\tilde{u}(x_{i_1},x_{i_2},x_{i_3}) - 3 \mu^2 \sin(\mu x_{i_1})\sin(\mu x_{i_2})\sin(\mu x_{i_3})|^2 \\
& + \frac{1}{M}\sum_{j=1}^M | \tilde{u}(x_{j_1},x_{j_2}, x_{j_3}) - \sin(\mu x_{j_1})\sin(\mu x_{j_2})\sin(\mu x_{j_3})|^2. \\
\end{align*}
\begin{figure}[H]
    \centering
    \subfloat[Loss During Fine-tuning]{%
        \includegraphics[width=6cm]{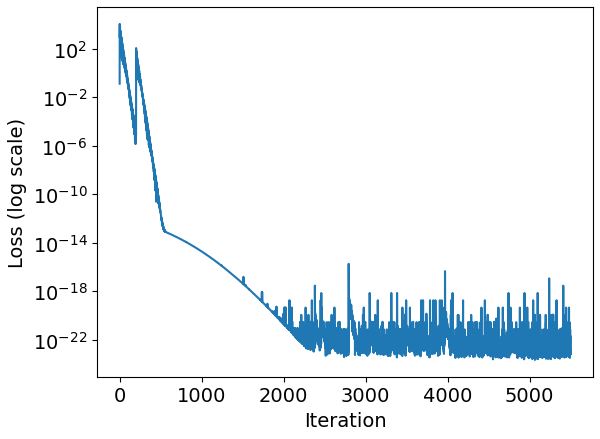}%
    }\qquad
    \subfloat[Relative Error During Fine-tuning]{%
        \includegraphics[width=6cm]{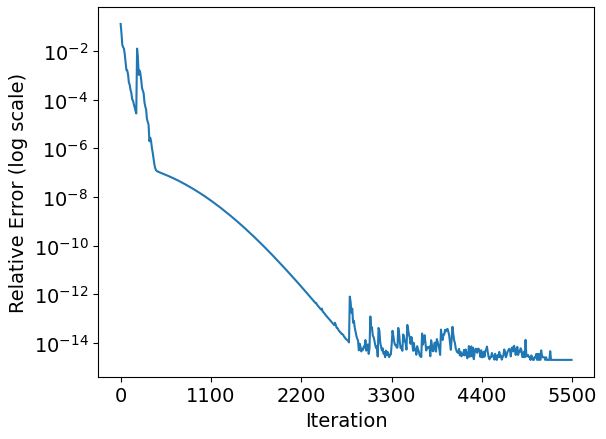}%
    }%
    \caption{Optimization profile for PDE~\eqref{eqn:Poisson6.1.4ex1}. 
    \textbf{(a)} Training loss of the candidate solution during fine-tuning. 
    \textbf{(b)} Absolute relative error of the candidate solution.}
    \label{fig:Poisson6.1.4ex1optim}
\end{figure}

The optimization profile indicates that FEX efficiently handles the complex geometry of the domain. In this case, 5,000 interior points and 5,000 boundary points (i.e., points on the surfaces of the spheres and the walls of the cube) were sampled uniformly to construct the loss functional at each iteration. Among the 5,000 boundary points, approximately half were drawn from the boundary of the cube, while the remaining half were evenly distributed among the surfaces of the spheres. Because the dimensionality of this problem is sufficiently low, the parameter grouping step of FEX-PG was omitted. This omission not only simplifies the procedure for low-dimensional problems but also isolates and validates the individual modifications introduced earlier in this section, confirming that these specific enhancements enable FEX to successfully solve the new equations. The average relative $L^2$ error achieved by FEX was $4.1 \times 10^{-14}$, approaching the double-precision machine epsilon. This result compares favorably with the findings of \cite{Ziqi_Liu_2020}, where the relative $L^2$ error was on the order of $10^{-2}$.
\\ 

\textbf{Example 2.} The other example is given as 
\begin{align}
    -\Delta u(\textbf{x}) = \mu^2 e^{\sin( \mu x_1)+\sin(\mu x_2)+\sin(\mu x_3)}\big(\cos^2( \mu x_1)+\cos^2(\mu x_2)+\cos^2(\mu x_3) - \sin( \mu x_1) - \sin(\mu x_2) -\sin(\mu x_3)\big),
    \label{eqn:Poisson6.1.4ex2}
\end{align}
where the boundary condition on the sides of the cube and surface of the spherical holes is chosen appropriately such that the true solution is
\begin{equation*}
    u(\textbf{x}) = e^{\sin(\mu x_1) + \sin(\mu x_2) + \sin(\mu x_3)}.
\end{equation*}
Once again $\mu = 7\pi$.  Our loss functional based on the least squares idea in \eqref{eqn:leastsquare} becomes
\begin{align*}
\mathcal{L}(\tilde{u}) & := \frac{1}{N}\sum_{i=1}^N| -\Delta\tilde{u}(x_{i_1},x_{i_2},x_{i_3})\\ 
& - \mu^2 e^{\sin( \mu x_{i_1})+\sin(\mu x_{i_2})+\sin(\mu x_{i_3})}\big(\cos^2( \mu x_{i_1})+\cos^2(\mu x_{i_2})+\cos^2(\mu x_{i_3}) - \sin( \mu x_{i_1}) - \sin(\mu x_{i_2}) -\sin(\mu x_{i_3})\big)|^2 \\
& + \frac{1}{M}\sum_{j=1}^M | \tilde{u}(x_{j_1},x_{j_2}, x_{j_3}) - e^{\sin(\mu x_{j_1})+\sin(\mu x_{j_2})+\sin(\mu x_{j_3})}|^2. \\
\end{align*}

This test examines how FEX-PG performs when the solution exhibits greater complexity, such as an exponential function applied to high-frequency sine components rather than their product. 

\begin{figure}[H]
    \centering
    \subfloat[Loss During Fine-tuning]{%
        \includegraphics[width=6cm]{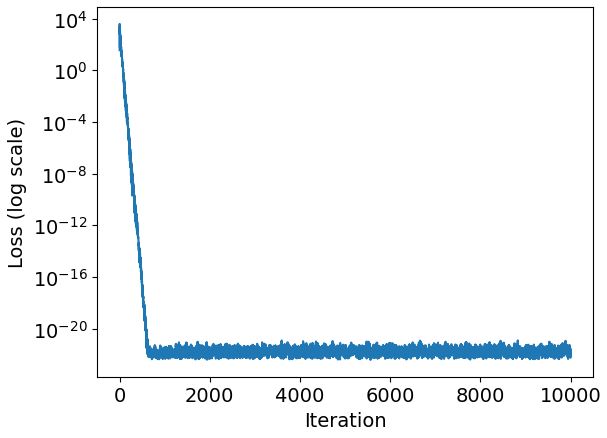}%
    }\qquad
    \subfloat[Relative Error During Fine-tuning]{%
        \includegraphics[width=6cm]{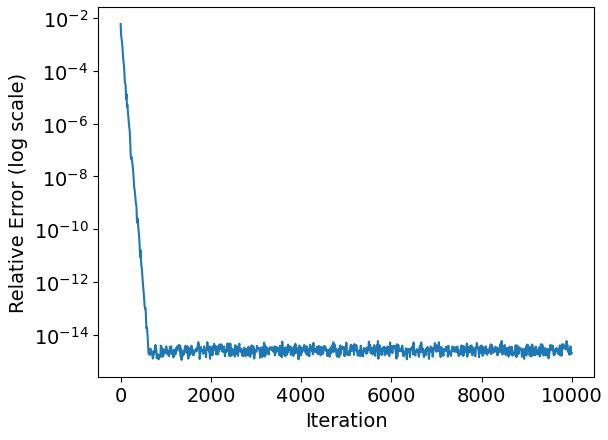}%
    }%
    \caption{Optimization profile for PDE~\eqref{eqn:Poisson6.1.4ex2} on the domain shown in Figure~\ref{fig:Poisson6.1.4domain}. 
    \textbf{(a)} Training loss of the candidate solution during fine-tuning. 
    \textbf{(b)} Absolute relative error of the candidate solution.}
    \label{fig:Poisson6.1.4ex2Optim}
\end{figure}

Again, we observe that the complex domain geometry poses no difficulty for problems of this type. Building on the success of the previous example, only 2,500 interior points and 2,500 boundary points were used per training iteration, enabling faster computation. To ensure adequate convergence, the number of iterations was doubled---although this proved unnecessary. The relative $L^2$ error averaged $3.2 \times 10^{-15}$, effectively reaching double-precision accuracy. This performance compares very favorably with the results reported in \cite{Zhi_Qin_John_Xu_2020}, where the relative $L^2$ error was on the order of $10^{-1}$. As shown in Figure~\ref{fig:Poisson6.1.4ex2Optim}, the optimization problem during fine-tuning is nearly trivial; at such low dimensionality, this is expected since the number of parameters is small. The real difficulty in problems of this kind lies primarily in the search phase---specifically, in solving the combinatorial optimization (CO) problem described in Section~\ref{sec:fex}.

\subsection{Eigenvalue Problem}
The last problem incorporates elements from many of the past examples, and adds the additional complexity of solving for the eigenvalue and function pair simultaneously.  Here we solve the Laplace eigenvalue problem with a zero boundary condition:
\begin{equation}
    \Delta u(\textbf{x}) = \lambda u(\textbf{x}) \text{, } u|_{\partial \Omega} = 0.
    \label{eqn:laplaceeigenvalue}
\end{equation}
Given the zero boundary condition, the equation admits a solution of the form
\begin{align}
    u(\textbf{x}) = \prod_{i=1}^d \sin(\frac{\pi x_i}{L})\text{, } \lambda = d\frac{\pi^2}{L}.
    \label{eqn:eigtruesol}
\end{align}
Importantly, the eigenfunction corresponding to any given eigenvalue is unique up to a constant multiple. In other words, if $u$ solves \eqref{eqn:laplaceeigenvalue} with eigenvalue $\lambda$ and $c \in \mathbb{R}$, then $cu$ is also a solution of \eqref{eqn:laplaceeigenvalue} associated with the same eigenvalue $\lambda$. It is worth noting that infinitely many eigenpairs satisfy \eqref{eqn:laplaceeigenvalue}; however, we focus on the eigenpair \eqref{eqn:eigtruesol}, which corresponds to the smallest eigenvalue. In practice, FEX often identifies the first two or three eigenpairs, though the smallest one is typically the final output, as its lower frequency makes the corresponding parameters easier to learn. To find the solution, we implement \eqref{eqn:eigenloss} and obtain the following loss functional:

\begin{equation*}
    \mathcal{L}(u) := \frac{1}{N}\sum_{i=1}^N|\Delta\tilde{u}(\textbf{x}_i) - \lambda \tilde{u}(\textbf{x}_i)|^2+ \alpha_b \frac{1}{M}\sum_{j=1}^M | \tilde{u}(\textbf{x}_j)|^2 + \alpha_n \text{min}_{i \in N}\{(|\tilde{u}(\textbf{x}_i)|^p -c)^2\},
\end{equation*}
We perform one hundred iterations of the operator-searching loop, conducting coarse tuning over the set of parameters that now includes $\lambda$, initialized as in \eqref{eqn:lambda_init}. The hyperparameters $\alpha_b$ and $\alpha_n$ in the loss functional are both set to 100. Following \cite{cai2023deepmartnetmartingalebased}, we solve this eigenvalue problem on a ten-dimensional cube.

\begin{figure}[H]
    \centering
    \subfloat[Loss During Fine-tuning]{%
        \includegraphics[width=6cm]{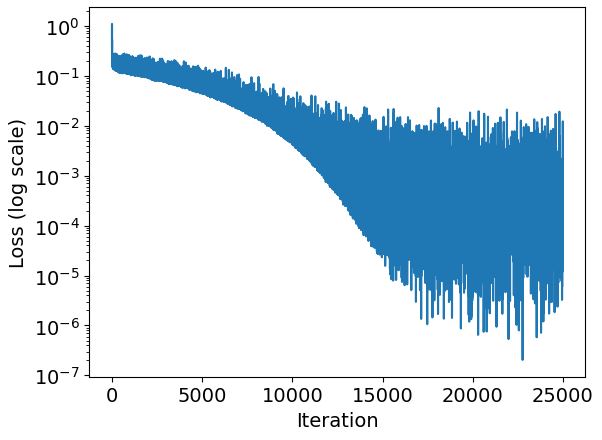}%
    }\qquad
    \subfloat[Relative Error During Fine-tuning]{%
        \includegraphics[width=6cm]{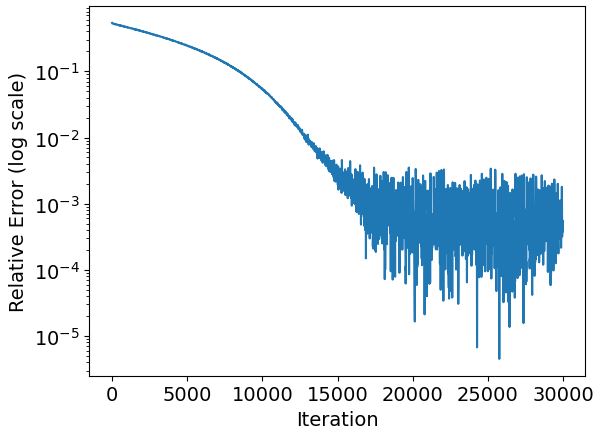}%
    }%
    \caption{Optimization profile for PDE~\eqref{eqn:laplaceeigenvalue} on the second domain. 
    \textbf{(a)} Training loss of the candidate solution during fine-tuning. 
    \textbf{(b)} Relative error of the candidate solution.}
    \label{fig:EigenvalueProblem}
\end{figure}

The noise observed in the optimization profiles in Figure~\ref{fig:EigenvalueProblem} arises primarily from the use of the minimum operator in \eqref{eqn:eigenloss}. Nonetheless, this formulation enables FEX-PG to consistently and efficiently identify a high-quality solution in every run. The average relative $L^2$ error of the resulting solution was $3 \times 10^{-3}$, representing an order of magnitude improvement in accuracy compared to \cite{cai2023deepmartnetmartingalebased}. As a preliminary proof of concept, this example demonstrates that the strengths of FEX-PG readily extend to eigenvalue problems. For reference, an example of the final solution produced by FEX-PG is shown below. Notably, because the frequency parameters within the sine function were grouped during the parameter-grouping step, the resulting function can be expressed exactly as a product across the dimensions of $\mathbf{x}$:

\begin{equation*}
    u(\textbf{x}) = \prod_{i=1}^{10} \sin(3.14168x_i)\text{, } \lambda = 98.69143.
\end{equation*}
The exact solution is
\begin{equation*}
    u(\textbf{x}) = \prod_{i=1}^{10} \sin(\pi x_i)\text{, } \lambda = 10\pi^2.
\end{equation*}
The accuracy achieved in this example is satisfactory, though not as high as in other cases where single- or double-precision accuracy was obtained. This reduction in accuracy is likely attributable to the inherent twofold optimization challenge of the eigenvalue problem: the parameters $\lambda$ and $u(\mathbf{x})$ are mutually dependent, and this coupling complicates the search for the global minimum of the loss functional.

\begin{figure}[H]
    \centering
    \includegraphics[width = 6cm]{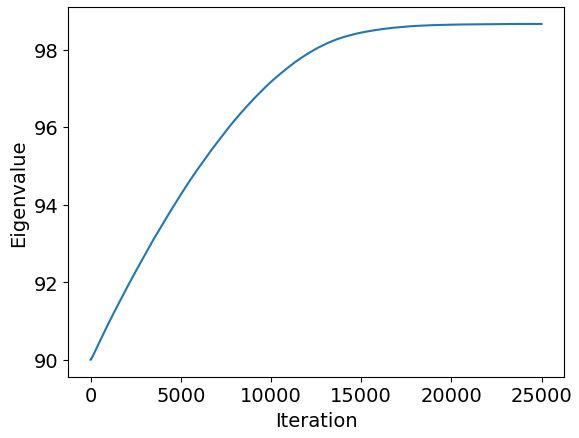}
    \caption{Convergence of the eigenvalue $\lambda$ over the iterations of fine-tuning}
    \label{fig:MartNet5.2.1eig}
\end{figure} 
As shown in Figure~\ref{fig:MartNet5.2.1eig}, the value of $\lambda$ is initialized at 90, following the estimate proposed in \eqref{eqn:lambda_init}. The raw candidate function provided by the controller was $u(\mathbf{x}) = \prod_{i=0}^{10}\sin(3x_i)$. These constitute excellent starting points when compared with the true solution \eqref{eqn:eigtruesol}, indicating that the controller was effectively trained through its searching loop. Nevertheless, it is evident that the eigenvalue problem presents a greater challenge than the previous examples, as reflected in the larger number of fine-tuning iterations required for convergence. Some improvement is likely achievable—for instance, employing distinct optimizers or learning rates for $\lambda$ and the parameters of $u(\mathbf{x})$ could enhance efficiency or accuracy—but such refinements are left for future work. \\

\subsection{Comparision Summary}
As a summary, the following Table \ref{tab:results_table} summarizes and compares all the numerical errors in the relative $L^2$ sense in the preceding tests in the numerical section to those in \cite{cai2023deepmartnetmartingalebased,Ziqi_Liu_2020}.

\begin{table}[H]
\begin{center}
\begin{tabular}{||c|| c c ||} 
 \hline
 Problem & FEX  & v.s. NN in \cite{cai2023deepmartnetmartingalebased,Ziqi_Liu_2020}\\ [0.5ex] 
 \hline\hline
 100-D Poisson-Boltzmann \eqref{eqn:PB_ex1} & ~10e-7 & ~5e-3 \cite{cai2023deepmartnetmartingalebased} \\ 
 \hline
 10-D Poisson-Boltzmann \eqref{eqn:PB_ex2} & 3.3e-6 & 2.5e-1 \cite{cai2023deepmartnetmartingalebased} \\
 \hline
 2-D Poisson, Small Hole Domain \eqref{eqn:Poisson6.1.3} & 4.9e-7 & 1e-2 \cite{Ziqi_Liu_2020} \\
 \hline
 2-D Poisson, Large Hole Domain\eqref{eqn:Poisson6.1.3} & 8.6e-7 &  8e-3 \cite{Ziqi_Liu_2020}\\
 \hline
 3-D Poisson, Multiplicative Solution \eqref{eqn:Poisson6.1.4ex1}& 4.1e-14 & 1e-2 \cite{Ziqi_Liu_2020}\\ [1ex] 
 \hline
 3-D Poisson, Exponential Solution \eqref{eqn:Poisson6.1.4ex2}& 3.2e-15 & 1e-0 \cite{Ziqi_Liu_2020} \\ [1ex] 
 \hline
 10-D Laplace Eigenvalue Problem \eqref{eqn:laplaceeigenvalue}& 3e-3 & 2.5e-1 \cite{cai2023deepmartnetmartingalebased} \\ [1ex] 
 \hline
\end{tabular}
\end{center}
\caption{A summery of comparisons from all numerical tests.}
\label{tab:results_table}
\end{table}

\subsection{Conclusion and Discussion}

This paper presents several novel developments that build upon the FEX-PG framework \cite{hardwick2024solvinghighdimensionalpartialintegral} and extend its applicability to a broader class of challenging problems. Specifically, we introduce significant enhancements to the Finite Expression Method (FEX) that enable it to address two notoriously difficult regimes in the numerical solution of partial differential equations (PDEs): highly oscillatory solutions and domains with complex geometries (e.g., regions with multiple holes). By incorporating a symbolic frequency composition module, a redesigned linear input layer, and new capabilities for solving eigenvalue problems, we substantially increase both the expressiveness and versatility of the FEX approach.

Across a diverse suite of benchmark problems—including nonlinear Poisson--Boltzmann equations, high-frequency Helmholtz-type problems on geometrically perforated domains, and high-dimensional eigenvalue problems—FEX demonstrates consistent accuracy and interpretability. The method achieves robust performance across varying dimensionalities and domain complexities, often producing errors several orders of magnitude smaller than state-of-the-art neural network solvers such as those reported in \cite{Ziqi_Liu_2020} and \cite{cai2023deepmartnetmartingalebased}.

A particularly notable advance is the introduction of the symbolic frequency composition module, which enables FEX to identify and combine the correct spectral components of a solution. In conjunction with the parameterized input layer, this module allows FEX to recover high-frequency solutions that are typically inaccessible to standard neural networks due to the F-Principle. For instance, FEX successfully recovers solutions such as $u(\textbf{x}) = \sin(7\pi x_0)\sin(7\pi x_1)$ (see \eqref{eqn:FEXoutputPoisson6.1.3}), with frequency estimates accurate to machine precision. Moreover, the symbolic representation of the resulting expressions, exemplified again in \eqref{eqn:FEXoutputPoisson6.1.3}, underscores one of FEX’s enduring strengths---interpretability---a property largely absent in black-box deep learning approaches.

In the context of eigenvalue problems, we demonstrated that FEX can recover both eigenfunctions and eigenvalues with strong accuracy by extending the loss functional and introducing a principled initialization scheme for $\lambda$. Although the achieved accuracy does not yet reach machine precision, the results are encouraging given the inherent difficulty of simultaneous optimization over the coupled eigenfunction and eigenvalue spaces.

Despite these advances, several open challenges remain. The expression-searching process (the reinforcement learning loop) remains computationally intensive, particularly for high-dimensional or stiff problems, and its performance depends strongly on the richness of the operator set. Problems characterized by widely separated frequency components---for example, PDEs with true solutions of the form $u(x_1,x_2,x_3) = \sin(10x_1)\sin(20x_2)\sin(30x_3)$---pose particular difficulties, as such frequency separation introduces instability in the reinforcement learning dynamics. Furthermore, scaling FEX to large-scale engineering PDEs, such as the Navier--Stokes or elastodynamic systems, will require further innovations in search efficiency and expression composition. Addressing these challenges represents a natural next step toward integrating FEX into mainstream computational science and engineering workflows.

\section*{Acknowledgement}
H. Y. was partially supported by the US National Science Foundation under awards DMS-2244988, the Office of Naval Research Award N00014-23-1-2007, and the DARPA D24AP00325-00. Approved for public release; distribution is unlimited. \\ \\

\bibliographystyle{plain} 
\bibliography{ref.bib}

\begin{thebibliography}{10}

\bibitem{cai2023deepmartnetmartingalebased}
Wei Cai, Andrew He, and Daniel Margolis.
\newblock Deepmartnet -- a martingale based deep neural network learning method
  for dirichlet bvps and eigenvalue problems of elliptic pdes in $r^d$, 2023.

\bibitem{cai2019phaseshiftdeepneural}
Wei Cai, Xiaoguang Li, and Lizuo Liu.
\newblock A phase shift deep neural network for high frequency approximation
  and wave problems, 2019.

\bibitem{chen2023deep}
Ke~Chen, Chunmei Wang, and Haizhao Yang.
\newblock Deep operator learning lessens the curse of dimensionality for
  {PDE}s.
\newblock {\em Transactions on Machine Learning Research}, 2023.

\bibitem{chizat}
Lenaic Chizat and Francis Bach.
\newblock Implicit bias of gradient descent for wide two-layer neural networks
  trained with the logistic loss, 2020.

\bibitem{2023JCoPh.47611862D}
Cuiyang {Ding}, Yijing {Zhou}, Wei {Cai}, Xuan {Zeng}, and Changhao {Yan}.
\newblock {A path integral Monte Carlo (PIMC) method based on Feynman-Kac
  formula for electrical impedance tomography}.
\newblock {\em Journal of Computational Physics}, 476:111862, March 2023.

\bibitem{E2017}
Weinan E, Jiequn Han, and Arnulf Jentzen.
\newblock Deep learning-based numerical methods for high-dimensional parabolic
  partial differential equations and backward stochastic differential
  equations.
\newblock {\em Communications in Mathematics and Statistics}, 5(4):349--380,
  Dec 2017.

\bibitem{e2017deepritzmethoddeep}
Weinan E and Bing Yu.
\newblock The deep ritz method: A deep learning-based numerical algorithm for
  solving variational problems, 2017.

\bibitem{fletcher2013practical}
Roger Fletcher.
\newblock {\em Practical methods of optimization}.
\newblock John Wiley \& Sons, 2013.

\bibitem{haberman}
Richard Haberman.
\newblock {\em Applied Partial Differential Equations With Fourier Series and
  Boundary Value Problems}.
\newblock Prentice Hall, Englewood Cliffs, NJ, 1983.

\bibitem{doi:10.1073/pnas.1718942115}
Jiequn Han, Arnulf Jentzen, and Weinan E.
\newblock Solving high-dimensional partial differential equations using deep
  learning.
\newblock {\em Proceedings of the National Academy of Sciences},
  115(34):8505--8510, 2018.

\bibitem{hardwick2024solvinghighdimensionalpartialintegral}
Gareth Hardwick, Senwei Liang, and Haizhao Yang.
\newblock Solving high-dimensional partial integral differential equations: The
  finite expression method, 2024.

\bibitem{doi:10.1142/S021953052350015X}
Yuling Jiao, Yanming Lai, Yisu Lo, Yang Wang, and Yunfei Yang.
\newblock Error analysis of deep ritz methods for elliptic equations.
\newblock {\em Analysis and Applications}, 22(01):57--87, 2024.

\bibitem{JiaoLaiWangYangYang2023-DGMW}
Yuling Jiao, Yanming Lai, Yang Wang, Haizhao Yang, and Yunfei Yang.
\newblock Convergence analysis of the deep galerkin method for weak solutions.
\newblock In Patricia Alonso~Ruiz, Michael Hinz, Kasso~A. Okoudjou, Luke~G.
  Rogers, and Alexander Teplyaev, editors, {\em From Classical Analysis to
  Analysis on Fractals: A Tribute to Robert Strichartz, Volume 1}, Applied and
  Numerical Harmonic Analysis, pages 53--82. Birkh{\"a}user, Cham, 2023.

\bibitem{Karniadakis2021}
George~Em Karniadakis, Ioannis~G. Kevrekidis, Lu~Lu, Paris Perdikaris, Sifan
  Wang, and Liu Yang.
\newblock Physics-informed machine learning.
\newblock {\em Nature Reviews Physics}, 3(6):422--440, Jun 2021.

\bibitem{KHOO_LU_YING_2021}
Yuehaw Khoo, Jianfeng Lu, and Lexing Ying.
\newblock Solving parametric pde problems with artificial neural networks.
\newblock {\em European Journal of Applied Mathematics}, 32(3):421–435, 2021.

\bibitem{Li_2022}
Haoya Li and Lexing Ying.
\newblock A semigroup method for high dimensional elliptic pdes and eigenvalue
  problems based on neural networks.
\newblock {\em Journal of Computational Physics}, 453:110939, March 2022.

\bibitem{liang2022stiffness}
Senwei Liang, Zhongzhan Huang, and Hong Zhang.
\newblock Stiffness-aware neural network for learning hamiltonian systems.
\newblock In {\em International Conference on Learning Representations}, 2022.

\bibitem{liang2024solving}
Senwei Liang, Shixiao~W Jiang, John Harlim, and Haizhao Yang.
\newblock Solving pdes on unknown manifolds with machine learning.
\newblock {\em Applied and Computational Harmonic Analysis}, 71:101652, 2024.

\bibitem{liang2024reproducing}
Senwei Liang, Liyao Lyu, Chunmei Wang, and Haizhao Yang.
\newblock Reproducing activation function for deep learning.
\newblock {\em Communications in Mathematical Sciences}, 22(2):285--314, 2024.

\bibitem{liang}
Senwei Liang and Haizhao Yang.
\newblock Finite expression method for solving high-dimensional partial
  differential equations.
\newblock {\em Journal of Machine Learning Research}, 26(138):1--31, 2025.

\bibitem{Ming_2021}
Yulei Liao and Pingbing Ming.
\newblock Deep nitsche method: Deep ritz method with essential boundary
  conditions.
\newblock {\em Communications in Computational Physics}, 29(5):1365–1384,
  June 2021.

\bibitem{JMLR:v25:22-0719}
Hao Liu, Haizhao Yang, Minshuo Chen, Tuo Zhao, and Wenjing Liao.
\newblock Deep nonparametric estimation of operators between infinite
  dimensional spaces.
\newblock {\em Journal of Machine Learning Research}, 25(24):1--67, 2024.

\bibitem{LiuHuangProtopapas2023_ResidualBasedErrorBound_PINNs}
Shuheng Liu, Xiyue Huang, and Pavlos Protopapas.
\newblock Residual-based error bound for physics-informed neural networks.
\newblock {\em arXiv preprint arXiv:2306.03786}, 2023.

\bibitem{Ziqi_Liu_2020}
Ziqi Liu, Wei Cai, and Zhi-Qin~John Xu.
\newblock Multi-scale deep neural network (mscalednn) for solving
  poisson-boltzmann equation in complex domains.
\newblock {\em Communications in Computational Physics}, 28(5):1970–2001,
  January 2020.

\bibitem{lu}
Liwei Lu, Hailong Guo, Xu~Yang, and Yi~Zhu.
\newblock Temporal difference learning for high-dimensional pides with jumps.
\newblock {\em SIAM Journal on Scientific Computing}, 46(4):C349--C368, 2024.

\bibitem{lu2021machinelearningellipticpdes}
Yiping Lu, Haoxuan Chen, Jianfeng Lu, Lexing Ying, and Jose Blanchet.
\newblock Machine learning for elliptic pdes: Fast rate generalization bound,
  neural scaling law and minimax optimality, 2021.

\bibitem{pmlr-v134-lu21a}
Yulong Lu, Jianfeng Lu, and Min Wang.
\newblock A priori generalization analysis of the deep ritz method for solving
  high dimensional elliptic partial differential equations.
\newblock In Mikhail Belkin and Samory Kpotufe, editors, {\em Proceedings of
  Thirty Fourth Conference on Learning Theory}, volume 134 of {\em Proceedings
  of Machine Learning Research}, pages 3196--3241. PMLR, 15--19 Aug 2021.

\bibitem{LuoYang2024-HNA}
Tao Luo and Haizhao Yang.
\newblock Two-layer neural networks for partial differential equations:
  optimization and generalization theory.
\newblock In Siddhartha Mishra and Alex Townsend, editors, {\em Numerical
  Analysis Meets Machine Learning}, volume~25 of {\em Handbook of Numerical
  Analysis}, pages 515--554. Elsevier, 2024.

\bibitem{maiti2025optimalneuralnetworkapproximation}
Ayan Maiti, Michelle Michelle, and Haizhao Yang.
\newblock Optimal neural network approximation for high-dimensional continuous
  functions.
\newblock {\em arxiv:2409.02363}, 2025.

\bibitem{MishraMolinaro2023_IMAJNA_forwardPINNs}
Siddhartha Mishra and Roberto Molinaro.
\newblock Estimates on the generalization error of physics-informed neural
  networks for approximating {PDE}s.
\newblock {\em IMA Journal of Numerical Analysis}, 43(1):1--43, 2023.

\bibitem{na2025cursedimensionalityneuralnetwork}
Sanghoon Na and Haizhao Yang.
\newblock Curse of dimensionality in neural network optimization.
\newblock {\em arxiv:2502.05360}, 2025.

\bibitem{petersen2021deep}
Brenden~K Petersen, Mikel~Landajuela Larma, Terrell~N. Mundhenk, Claudio~Prata
  Santiago, Soo~Kyung Kim, and Joanne~Taery Kim.
\newblock Deep symbolic regression: Recovering mathematical expressions from
  data via risk-seeking policy gradients.
\newblock In {\em International Conference on Learning Representations}, 2021.

\bibitem{pmlr-v97-rahaman19a}
Nasim Rahaman, Aristide Baratin, Devansh Arpit, Felix Draxler, Min Lin, Fred
  Hamprecht, Yoshua Bengio, and Aaron Courville.
\newblock On the spectral bias of neural networks.
\newblock In Kamalika Chaudhuri and Ruslan Salakhutdinov, editors, {\em
  Proceedings of the 36th International Conference on Machine Learning},
  volume~97 of {\em Proceedings of Machine Learning Research}, pages
  5301--5310. PMLR, 09--15 Jun 2019.

\bibitem{raissi}
M.~Raissi, P.~Perdikaris, and G.E. Karniadakis.
\newblock Physics-informed neural networks: A deep learning framework for
  solving forward and inverse problems involving nonlinear partial differential
  equations.
\newblock {\em Journal of Computational Physics}, 378:686--707, 2019.

\bibitem{Shen_2021}
Zuowei Shen, Haizhao Yang, and Shijun Zhang.
\newblock Deep network with approximation error being reciprocal of width to
  power of square root of depth.
\newblock {\em Neural Computation}, 33(4):1005–1036, 2021.

\bibitem{Shen_2021_3layers}
Zuowei Shen, Haizhao Yang, and Shijun Zhang.
\newblock Neural network approximation: Three hidden layers are enough.
\newblock {\em Neural Networks}, 141:160–173, September 2021.

\bibitem{Sirignano_2018}
Justin Sirignano and Konstantinos Spiliopoulos.
\newblock Dgm: A deep learning algorithm for solving partial differential
  equations.
\newblock {\em Journal of Computational Physics}, 375:1339–1364, December
  2018.

\bibitem{song}
Zezheng Song, Maria~K. Cameron, and Haizhao Yang.
\newblock A finite expression method for solving high-dimensional committor
  problems.
\newblock {\em SIAM Journal of Scientific Computing}, 2024.

\bibitem{song2024finiteexpressionmethodlearning}
Zezheng Song, Chunmei Wang, and Haizhao Yang.
\newblock Finite expression method for learning dynamics on complex networks.
\newblock {\em arxiv:2401.03092}, 2024.

\bibitem{9321497}
Stephan Wojtowytsch and Weinan E.
\newblock Can shallow neural networks beat the curse of dimensionality? a mean
  field training perspective.
\newblock {\em IEEE Transactions on Artificial Intelligence}, 1(2):121--129,
  2020.

\bibitem{Zhi_Qin_John_Xu_2020}
Zhi-Qin~John Xu, Yaoyu Zhang, Tao Luo, Yanyang Xiao, and Zheng Ma.
\newblock Frequency principle: Fourier analysis sheds light on deep neural
  networks.
\newblock {\em Communications in Computational Physics}, 28(5):1746–1767,
  January 2020.

\bibitem{shen2022deepnetworkapproximationachieving}
Shijun Zhang, Zuowei Shen, and Haizhao Yang.
\newblock Deep network approximation: Achieving arbitrary accuracy with fixed
  number of neurons.
\newblock {\em Journal of Machine Learning Research}, 23(276):1--60, 2022.

\bibitem{zhang2021fbsdebasedneuralnetwork}
Wenzhong Zhang and Wei Cai.
\newblock Fbsde based neural network algorithms for high-dimensional
  quasilinear parabolic pdes, 2021.

\end{thebibliography}
\end{document}